\documentclass[12pt]{article}

\usepackage{amsfonts}
\usepackage{amsmath}
\usepackage{amssymb,euscript}

\newcommand{\Xcomment}[1]{}

\makeatletter
\renewcommand{\section}{\@startsection{section}{1}{0pt}%
{-3.5ex plus -1ex minus -.2ex}{2.3ex plus .2ex}%
{\normalfont\Large}}
\makeatother

 \makeatletter
\renewcommand{\subsection}{\@startsection{subsection}{2}{0pt}%
{-3.0ex plus -1ex minus -.2ex}{1.5ex plus .2ex}%
{\normalfont\normalsize\bf}}
 \makeatother

 \makeatletter
\renewcommand{\subsubsection}{\@startsection{subsubsection}{3}{0pt}%
{-3.0ex plus -1ex minus -.2ex}{-1.5ex plus -.2ex}%
{\normalfont\normalsize\bf}}
 \makeatother

\newtheorem{theorem}{Theorem}[section]
\newtheorem{lemma}[theorem]{Lemma}

\newtheorem{prop}[theorem]{Proposition}

{$\qed$\medskip}

\makeatletter \@addtoreset{equation}{section} \makeatother

\newenvironment{numitem1}{\refstepcounter{equation}\begin{enumerate}%
\item[(\thesection.\arabic{equation})]}{\end{enumerate}}

\newcommand{\refeq}[1]{(\ref{eq:#1})}  

\def\rest#1{_{\,\vrule height 1.8ex width 0.05em depth 0pt\, #1}}
\def\qed{ \ \vrule width.1cm height.3cm depth0cm}

\def\tilde{\widetilde}
\def\hat{\widehat}
\def\eps{\varepsilon}

\def\Rset{{\mathbb R}}
\def\Zset{{\mathbb Z}}

\def\Bscr{\mathcal{B}}

\def\Escr{\mathcal{E}}
\def\Fscr{\mathcal{F}}
\def\Gscr{\mathcal{G}}

\def\Kscr{\mathcal{K}}

\def\Mscr{\mathcal{M}}

\def\Kmn{\mathcal{K}^{(-n)}}
\def\Kmone{\mathcal{K}^{(-1)}}

\def\prv{\check{v}}

\def\Kup{K^{\uparrow}}
\def\Klow{K^{\downarrow}}
\def\Piup{\Pi^{\uparrow}}
\def\Pilow{\Pi^{\downarrow}}
\def\vup{v^{\uparrow}}
\def\vlow{v^{\downarrow}}

\def\parup{c^{\uparrow}}
\def\parlow{c^{\downarrow}}
\def\heartup{\hslash^{\uparrow}}
\def\heartlow{\hslash^{\downarrow}}
\def\parmid{c^{\uparrow\downarrow}}

\def\Kmid{K^{\uparrow\downarrow}}
\def\heartmid{\hslash^{\uparrow\downarrow}}
\def\bmid{\hslash^{\downarrow\uparrow}}

\def\bfa{a}
\def\bfb{b}
\def\bfc{c}

\def\bfp{p}
\def\bfq{q}

\def\bfzero{{\bf 0}}

\def\Pin{P^{\rm in}}
\def\Pout{P^{\rm out}}
\def\ellin{t}
\def\ellout{h}

\def\eNW{e^{\rm NW}}
\def\eSW{e^{\rm SW}}
\def\eNE{e^{\rm NE}}
\def\eSE{e^{\rm SE}}

\begin{document}

\begin{center}
 {\large\bf  Assembling crystals of type A}%
\footnote{Supported by RFBR grant 10-01-9311-CNRSL\_\,a.}
 \end{center}
 \medskip

\begin{center}
{\sc Vladimir~I.~Danilov}\footnote{Central Institute of Economics and
Mathematics of the RAS (47, Nakhimovskii Prospect, 117418 Moscow, Russia);
email: danilov@cemi.rssi.ru.},
{\sc Alexander~V.~Karzanov}\footnote{Institute for System Analysis of the RAS
(9, Prospect 60 Let Oktyabrya, 117312 Moscow, Russia); email:
sasha@cs.isa.ru.}\footnote{A part of this research was done while this author
was visiting Equipe Combinatoire et Optimisation, Univ. Paris-6, and Institut
f\"ur Diskrete Mathematik, Univ. Bonn. Corresponding author.}, \\
{\sc and Gleb~A.~Koshevoy}\footnote{Central Institute of Economics and
Mathematics of the RAS (47, Nakhimovskii Prospect, 117418 Moscow, Russia) and
Laboratoire J.-V.Poncelet (11, Bolshoy Vlasyevskiy Pereulok, 119002 Moscow,
Russia); email: koshevoy@cemi.rssi.ru.}
\end{center}

 \bigskip
 \begin{quote}
 {\bf Abstract.}
 {\small
Regular $A_n$-\emph{crystals} are certain edge-colored directed graphs which
are related to representations of the quantized universal enveloping algebra
$U_q(\mathfrak{sl}_{n+1})$. For such a crystal $K$ with colors $1,2,\ldots,n$,
we consider its maximal connected subcrystals with colors $1,\ldots,n-1$ and
with colors $2,\ldots,n$ and characterize the interlacing structure for all
pairs of these subcrystals. This is used to give a recursive description of the
combinatorial structure of $K$ and develop an efficient procedure of assembling
$K$.
 \smallskip

\smallskip

  }
  \end{quote}

\section{Introduction} \label{sec:intr}

{\em Crystals} are certain ``exotic'' edge-colored graphs. This graph-theoretic
abstraction, introduced by Kashiwara~\cite{kas-90,kas-95}, has proved its
usefulness in the theory of representations of Lie algebras and their quantum
analogues. In general, a finite \emph{crystal} is a finite directed graph $K$
such that: the edges are partitioned into $n$ subsets, or \emph{color} classes,
labeled $1,\ldots,n$, each connected monochromatic subgraph of $K$ is a simple
directed path, and there is a certain interrelation between the lengths of such
paths, which depends on the $n\times n$ Cartan matrix $M=(m_{ij})$ related to a
given Lie algebra $\mathfrak{g}$. Of most interest are crystals of
representations, or {\em regular} crystals. They are associated to elements of
a certain basis of the highest weight integrable modules (representations) over
the quantized universal enveloping algebra $U_q(\mathfrak{g})$.

This paper continues our combinatorial study of crystals begun
in~\cite{A2,cross} and considers $n$-colored regular crystals of type A, where
the number $n$ of colors is arbitrary. Recall that {\em type~A} concerns
$\mathfrak{g}=\mathfrak{sl}_{n+1}$; in this case the Cartan matrix $M$ is
viewed as $m_{ij}=-1$ if $|i-j|=1$, $m_{ij}=0$ if $|i-j|>1$, and $m_{ii}=2$. We
will refer to a regular $n$-colored crystal of type A as an $A_n$-{\em crystal}
and omit the index $n$ when the number of colors is not specified. Since we are
going to deal with finite regular crystals only, the adjectives ``finite'' and
``regular'' will usually be omitted. Also we assume that any crystal in
question is (weakly) \emph{connected}, i.e. it is not the disjoint union of two
nonempty graphs (which does not lead to loss of generality).

It is known that any A-crystal $K$ possesses the following properties. (i) $K$
is acyclic (i.e. has no directed cycles) and has exactly one zero-indegree
vertex, called the {\em source}, and exactly one zero-outdegree vertex, called
the {\em sink} of $K$. (ii) For any $I\subseteq\{1,\ldots,n\}$, each
(inclusion-wise) maximal connected subgraph of $K$ whose edges have colors from
$I$ is a crystal related to the corresponding $I\times I$ submatrix of the
Cartan matrix of $K$. Throughout, speaking of a \emph{subcrystal} of $K$, we
will always mean a subgraph of this sort.

Two-colored subcrystals are of especial importance, due to the result
in~\cite{KKM-92} that for a crystal (of any type) with exactly one
zero-indegree vertex, the regularity of all two-colored subcrystals implies the
regularity of the whole crystal. Let $K'$ be a two-colored subcrystal with
colors $i,j$ in an A-crystal $K$. Then $K'$ is the Cartesian product of a path
with color $i$ and a path with color $j$ (forming an $A_1\times A_1$-crystal)
when $|i-j|>1$, and an $A_2$-crystal when $|i-j|=1$. (The A-crystals belong to
the group of {\em simply-laced} crystals, which are characterized by the
property that each two-colored subcrystal is of type $A_1\times A_1$ or $A_2$;
for more details, see~\cite{Stem}.)

Another important fact is that for any $n$-tuple $c=(c_1,\ldots,c_n)$ of
nonnegative integers, there exists exactly one $A_n$-crystal $K$ such that each
$c_i$ is equal to the length of the maximal path with color $i$ beginning at
the source (for a short combinatorial proof, see~\cite[Sec.~2]{cross}). We
denote the crystal $K$ determined by $c$ in this way by $K(c)$, and refer to
$c$ as the {\em parameter} of this crystal.

There have been known several ways to define A-crystals; in particular, via
Gelfand--Tsetlin pattern model, semistandard Young tableaux, Littelmann's path
model; see~\cite{GT-50,KN-94,Lit-95,Littl}. In the last decade there appeared
additional, more enlightening, descriptions. A short list of ``local'' defining
axioms for A-crystals is pointed out in~\cite{Stem}, and an explicit
construction for $A_2$-crystals is given in~\cite{A2}. According to that
construction, any $A_2$-crystal can be obtained from an $A_1\times A_1$-crystal
by replacing each monochromatic path of the latter by a graph viewed as a
triangular half of a directed square grid.

When $n>2$, the combinatorial structure of $A_n$-crystals becomes much more
complicated, even for $n=3$. Attempting to learn more about this structure, we
elaborated in~\cite{cross} a new combinatorial construction, the so-called {\em
crossing model} (which is a refinement of the Gelfand--Tsetlin pattern model).
This powerful tool has helped us to reveal more structural features of an
$A_n$-crystal $K=K(c)$. In particular, $K$ has the so-called {\em principal
lattice}, a vertex subset $\Pi$ with the following nice properties:

(P1) $\Pi$ contains the source and sink of $K$, and the vertices $v\in \Pi$ are
bijective to the elements of the integer box $\Bscr(c):=\{a\in\Zset^n\colon
0\le a\le c\}$; we write $v=\prv[a]$;

(P2) For any $a,a'\in\Bscr(c)$ with $a\le a'$, the \emph{interval} of $K$ from
$\prv[a]$ to $\prv[a']$ (i.e. the subgraph of $K$ formed by the vertices and
edges contained in directed paths from $\prv[a]$ to $\prv[a']$) is isomorphic
to the $A_n$-crystal $K(a'-a)$, and its principal lattice consists of the
principal vertices $\prv[a'']$ of $K$ with $a\le a''\le a'$;

(P3) The set $\Kscr^{(-n)}$ of $(n-1)$-colored subcrystals $K'$ of $K$ having
colors $1,\ldots,n-1$ is bijective to $\Pi$; more precisely, $K'\cap\Pi$
consists of a single vertex, called the \emph{heart} of $K'$ w.r.t. $K$;
similarly for the set $\Kscr^{(-1)}$ of subcrystals of $K$ with colors
$2,\ldots,n$.
\smallskip

Note that a sort of ``principal lattice'' satisfying (P1) and (P2) can be
introduced for crystals of types B and C as well (see~\cite[Sec.~8]{ABC}), and
probably for the other classical types; however, (P3) does not remain true in
general for those types. Property~(P3) is crucial in our study of A-crystals in
this paper.

For $a\in\Bscr(c)$, let $\Kup[a]$ (resp. $\Klow[a]$) denote the subcrystal in
$\Kscr^{(-n)}$ (resp. in $\Kscr^{(-1)}$) that contains the principal vertex
$\prv[a]$; we call it the {\em upper} (resp. {\em lower}) subcrystal at $a$. It
is shown in~\cite{cross} that the parameter of this subcrystal is expressed by
a linear function of $c$ and $a$, and that the total amount of upper (lower)
subcrystals with a fixed parameter $c'$ is expressed by a piece-wise linear
function of $c$ and $c'$.
\smallskip

In this paper, we further utilize the crossing model, aiming to obtain a
refined description of the structure of an $A_n$-crystal $K$. We study the
intersections of subcrystals $\Kup[a]$ and $\Klow[b]$ for any $a,b\in\Bscr(c)$.
This intersection may be empty or consist of one or more subcrystals with
colors $2,\ldots,n-1$, called \emph{middle} subcrystals of $K$. Each of these
middle subcrystals $\tilde K$ is therefore a lower subcrystal of $\Kup[a]$ and
an upper subcrystal of $\Klow[b]$; so $\tilde K$ has a unique vertex $z$ in the
principal lattice $\Piup$ of the former, and a unique vertex $z'$ in the
principal lattice $\Pilow$ of the latter. Our main structural results --
Theorems~\ref{tm:two_div} and~\ref{tm:two_loci} -- give explicit relations
involving the ``loci'' $a$ and $b$ in $\Pi$, the ``deviation'' from $z$ of (the
heart of) $\tilde K$ in $\Piup$, and the ``deviation'' from $z'$ of $\tilde K$
in $\Pilow$.

This gives rise to a recursive procedure of assembling of the $A_n$-crystal
$K(c)$. More precisely, suppose that the $(n-1)$-colored crystals $\Kup[a]$ and
$\Klow[b]$ for all $a,b\in\Bscr(c)$ are already constructed. Then we can
combine these subcrystals to obtain the desired crystal $K(c)$, by properly
identifying the corresponding middle subcrystals (if any) for each pair
$\Kup[a],\Klow[b]$. This recursive method is implemented as an efficient
algorithm which, given a parameter $c\in\Zset_+^n$, outputs the crystal $K(c)$.
The running time of the algorithm and the needed space are bounded by
$Cn^2|K(c)|$, where $C$ is a constant and $|K(c)|$ is the size of $K(c)$. (It
may be of practical use for small $n$ and $c$; in general, an $A_n$-crystal has
``dimension'' $\frac{n(n+1)}{2}$ and its size grows sharply by increasing $c$.)

This paper is organized as follows. Section~\ref{sec:prelim} contains basic
definitions and backgrounds. Here we recall ``local'' axioms and the crossing
model for A-crystals, and review needed results on the principal lattice $\Pi$
of an $A_n$-crystal and relations between $\Pi$ and the $(n-1)$-colored
subcrystals from~\cite{cross}. Section~\ref{sec:ass_A} states
Theorems~\ref{tm:two_div} and~\ref{tm:two_loci} and discusses a recursive
description of the structure of an $A_n$-crystal $K$ and the algorithm of
assembling $K$. These theorems are proved in Section~\ref{sec:proof}. Our
assembling method is illustrated in Section~\ref{sec:illustr} with two special
cases of A-crystals: for an arbitrary $A_2$-crystal (in which case the method
can be compared with the explicit combinatorial construction in~\cite{A2}), and
for the particular $A_3$-crystal $K(1,1,1)$. \smallskip

It should be noted that the obtained structural results on A-crystals can also
be used to give a direct combinatorial proof of the known fact that any regular
$B_n$-crystal ($C_n$-crystal) can be extracted, in a certain way, from a
symmetric $A_{2n-1}$-crystal (resp. $A_{2n}$-crystal); this relationship is
discussed in detail in~\cite[Sec.~5--8]{ABC}. Here an $A_k$-crystal with
parameter $(c_1,\ldots,c_k)$ is called symmetric if $c_i=c_{k+1-i}$.

\section{Preliminaries} \label{sec:prelim}

In this section we recall ``local'' axioms defining (regular)  A-crystals,
explain the construction of crossing model, and review facts about the
principal lattice and subcrystals established in~\cite{cross} that will be
needed later.

 \subsection{A-crystals} \label{ssec:typeA}

Stembridge~\cite{Stem} pointed out a list of ``local'' graph-theoretic axioms
for the regular simply-laced crystals. The A-crystals form a subclass of those
and are defined by axioms (A1)--(A5) below; these axioms are given in a
slightly different, but equivalent, form compared with~\cite{Stem}.

Let $K=(V(K),E(K))$ be a directed graph whose edge set is partitioned into $n$
subsets $E_1,\ldots,E_n$, denoted as $K=(V(K),E_1\sqcup\ldots\sqcup E_n)$. We
assume $K$ to be (weakly) connected. We say that an edge $e\in E_i$ has {\em
color} $i$, or is an $i$-{\em edge}.

Unless explicitly stated otherwise, by a {\em path} we mean a simple
\emph{finite} directed path, i.e. a sequence of the form
$(v_0,e_1,v_1,\ldots,e_k,v_k)$, where $v_0,v_1,\ldots,v_k$ are distinct
vertices and each $e_i$ is an edge from $v_{i-1}$ to $v_i$ (admitting $k=0$).

The first axiom concerns the
structure of monochromatic subgraphs of $K$.
\begin{itemize}
\item[(A1)] For $i=1,\ldots,n$, each connected subgraph of $(V(K),E_i)$ is
a path.
  \end{itemize}

So each vertex of $K$ has at most one incoming $i$-edge and at most one
outgoing $i$-edge, and therefore one can associate to the set $E_i$ a {\em
partial invertible operator} $F_i$ acting on vertices: $(u,v)$ is an $i$-edge
if and only if $F_i$ {\em acts} at $u$ and $F_i(u)=v$ (or $u=F_i^{-1}(v)$,
where $F_i^{-1}$ is the partial operator inverse to $F_i$). Since $K$ is
connected, one can use the operator notation to express any vertex via another
one. For example, the expression $F_1^{-1}F_3^2F_2(v)$ determines the vertex
$w$ obtained from a vertex $v$ by traversing 2-edge $(v,v')$, followed by
traversing 3-edges $(v',u)$ and $(u,u')$, followed by traversing 1-edge
$(w,u')$ in backward direction. Emphasize that every time we use such an
operator expression in what follows, this automatically says that all
corresponding edges do exist in $K$.

We refer to a monochromatic path with color $i$ on the edges as an $i$-{\em
path}. So each maximal $i$-path is an $A_1$-subcrystal with color $i$ in $K$.
The maximal $i$-path passing a given vertex $v$ (possibly consisting of the
only vertex $v$) is denoted by $P_i(v)$, its part from the first vertex to $v$
by $\Pin_i(v)$, and its part from $v$ to the last vertex by $\Pout_i(v)$ (the
\emph{tail} and \emph{head} parts of $P$ w.r.t. $v$). The lengths (i.e. the
numbers of edges) of $\Pin_i(v)$ and $\Pout_i(v)$ are denoted by $\ellin_i(v)$
and $\ellout_i(v)$, respectively.

Axioms (A2)--(A5) concern interrelations of different colors $i,j$. They say
that each component of the two-colored graph $(V(K),E_i\sqcup E_j)$ forms an
$A_2$-crystal when colors $i,j$ are {\em neighboring}, which means that
$|i-j|=1$, and forms an $A_1\times A_1$-crystal otherwise.

When we traverse an edge of color $i$, the head and tail part lengths of
maximal paths of another color $j$ behave as follows:
 \begin{itemize}
 \item[(A2)]
For different colors $i,j$ and for an edge $(u,v)$ with color $i$, one
holds $\ellin_j(v)\le\ellin_j(u)$ and $\ellout_j(v)\ge\ellout_j(u)$. The
value $(\ellout_j(u)-\ellin_j(u))-(\ellout_j(v)-\ellin_j(v))$ is the
constant $m_{ij}$ equal to $-1$ if $|i-j|=1$, and 0 otherwise.
Furthermore, $h_j$ is convex on each $i$-path, in the sense that if
$(u,v),(v,w)$ are consecutive $i$-edges, then $h_j(u)+h_j(w)\ge 2h_j(v)$.
  \end{itemize}
These constants $m_{ij}$ are just the off-diagonal entries of the Cartan
$n\times n$ matrix $M$ related to the crystal type A and the number $n$ of
colors.

It follows that for neighboring colors $i,j$, each maximal $i$-path $P$
contains a unique vertex $r$ such that: when traversing any edge $e$ of $P$
before $r$ (i.e. $e\in\Pin_i(r)$), the tail length $\ellin_j$ decreases by 1
while the head length $\ellout_j$ does not change, and when traversing any edge
of $P$ after $r$, $\ellin_j$ does not change while $\ellout_j$ increases by 1.
This $r$ is called the {\em critical} vertex for $P,i,j$. To each $i$-edge
$e=(u,v)$ we associate {\em label} $\ell_{j}(e):=\ellout_j(v)-\ellout_j(u)$;
then $\ell_j(e)\in\{0,1\}$ and $t_j(v)=t_j(u)-1+\ell_j(e)$. Emphasize that the
critical vertices on a maximal $i$-path $P$ w.r.t. its neighboring colors
$j=i-1$ and $j=i+1$ may be different (and so are the edge labels on $P$).

Two operators $F=F_i^{\alpha}$ and $F'=F_j^\beta$, where $\alpha,\beta\in
\{1,-1\}$, are said to {\em commute} at a vertex $v$ if each of $F,F'$ acts at
$v$ (i.e. corresponding $i$-edge and $j$-edge incident with $v$ exist) and
$FF'(v)=F'F(v)$. The third axiom indicates situations when such operators
commute for neighboring $i,j$.

\begin{itemize}
\item[(A3)] Let $|i-j|=1$. (a) If a vertex $u$ has outgoing $i$-edge $(u,v)$
and outgoing $j$-edge $(u,v')$ and if $\ell_{j}(u,v)=0$, then
$\ell_{i}(u,v')=1$ and $F_i,F_j$ commute at $v$.
Symmetrically: (b) if a vertex $v$ has incoming $i$-edge $(u,v)$ and
incoming $j$-edge $(u',v)$ and if $\ell_{j}(u,v)=1$, then
$\ell_{i}(u',v)=0$ and $F_i^{-1},F_j^{-1}$ commute at $v$.
(See the picture.)
  \end{itemize}
 \begin{center}
  \unitlength=1mm
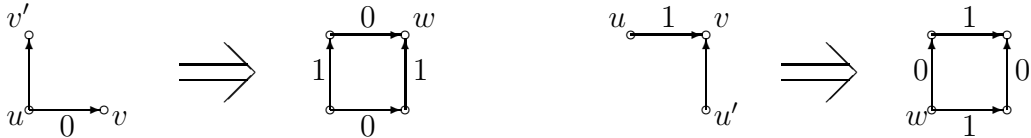
  \begin{picture}(140,20)
\put(5,5){\circle{1.0}}
\put(15,5){\circle{1.0}}
\put(45,5){\circle{1.0}}
\put(55,5){\circle{1.0}}
\put(95,5){\circle{1.0}}
\put(125,5){\circle{1.0}}
\put(135,5){\circle{1.0}}
\put(5,15){\circle{1.0}}
\put(45,15){\circle{1.0}}
\put(55,15){\circle{1.0}}
\put(85,15){\circle{1.0}}
\put(95,15){\circle{1.0}}
\put(125,15){\circle{1.0}}
\put(135,15){\circle{1.0}}
\put(5,5){\vector(1,0){9.5}}
\put(45,5){\vector(1,0){9.5}}
\put(125,5){\vector(1,0){9.5}}
\put(45,15){\vector(1,0){9.5}}
\put(85,15){\vector(1,0){9.5}}
\put(125,15){\vector(1,0){9.5}}
\put(5,5){\vector(0,1){9.5}}
\put(45,5){\vector(0,1){9.5}}
\put(55,5){\vector(0,1){9.5}}
\put(95,5){\vector(0,1){9.5}}
\put(125,5){\vector(0,1){9.5}}
\put(135,5){\vector(0,1){9.5}}
\put(25,9){\line(1,0){9}}
\put(25,11){\line(1,0){9}}
\put(105,9){\line(1,0){9}}
\put(105,11){\line(1,0){9}}
\put(31,6){\line(1,1){4}}
\put(31,14){\line(1,-1){4}}
\put(111,6){\line(1,1){4}}
\put(111,14){\line(1,-1){4}}
\put(2,3){$u$}
\put(2,16){$v'$}
\put(16,3){$v$}
\put(56,16){$w$}
\put(82,16){$u$}
\put(96,16){$v$}
\put(96,3){$u'$}
\put(121.5,3){$w$}
\put(9,1.5){0}
\put(49,1.5){0}
\put(129,1.5){1}
\put(49,16){0}
\put(89,16){1}
\put(129,16){1}
\put(42.5,9){1}
\put(56,9){1}
\put(122.5,9){0}
\put(136,9){0}
  \end{picture}
 \end{center}

One easily shows that if four vertices are connected by two $i$-edges $e,e'$
and two $j$-edges $\tilde e,\tilde e'$ (forming a ``square''), then
$\ell_{j}(e)=\ell_{j}(e')\ne \ell_{i}(\tilde e)=\ell_{i}(\tilde e')$ (as
illustrated in the picture). Another important consequence of (A3) is that for
neighboring colors $i,j$, if $v$ is the critical vertex on a maximal $i$-path
w.r.t. color $j$, then $v$ is also the critical vertex on the maximal $j$-path
passing $v$ w.r.t. color $i$, i.e. we can speak of common critical vertices for
the pair $\{i,j\}$.

The fourth axiom points out situations when, for neighboring $i,j$, the
operators $F_i,F_j$ and their inverse ones  ``remotely commute'' (forming the
``Verma relation of degree 4'').

\begin{itemize}
\item[(A4)] Let $|i-j|=1$.
(i) If a vertex $u$ has outgoing edges with color $i$ and color $j$ and if each
edge is labeled 1 w.r.t. the other color, then $F_iF_j^2F_i(u)=F_jF_i^2F_j(u)$.
Symmetrically: (ii) if $v$ has incoming edges with color $i$ and color $j$ and
if both are labeled 0, then $F_i^{-1}(F_j^{-1})^2
F_i^{-1}(v)=F_j^{-1}(F_i^{-1})^2 F_j^{-1}(v)$. (See the picture.)
  \end{itemize}
 \begin{center}
  \unitlength=1mm
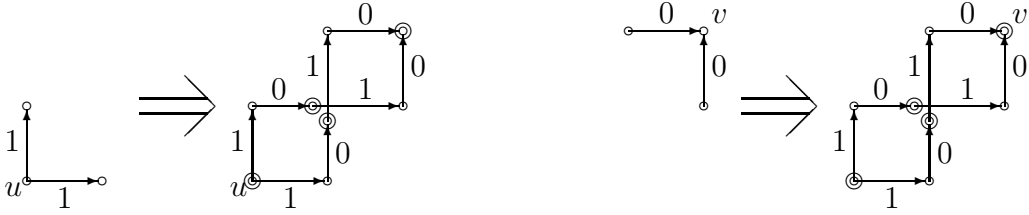
  \begin{picture}(147,30)
\put(5,5){\circle{1.0}}
\put(15,5){\circle{1.0}}
\put(5,15){\circle{1.0}}
\put(2,3){$u$}
\put(5,5){\vector(1,0){9.5}}
\put(5,5){\vector(0,1){9.5}}
\put(9,1.5){1}
\put(2,9){1}
\put(20,14){\line(1,0){9}}
\put(20,16){\line(1,0){9}}
\put(26,11){\line(1,1){4}}
\put(26,19){\line(1,-1){4}}
 \put(35,5){\circle{1.0}}
 \put(45,5){\circle{1.0}}
 \put(35,15){\circle{1.0}}
 \put(45,13){\circle{1.0}}
 \put(43,15){\circle{1.0}}
 \put(55,15){\circle{1.0}}
 \put(45,25){\circle{1.0}}
 \put(55,25){\circle{1.0}}
 \put(35,5){\circle{2.0}}
 \put(45,13){\circle{2.0}}
 \put(43,15){\circle{2.0}}
 \put(55,25){\circle{2.0}}
 \put(35,5){\vector(1,0){9.5}}
 \put(35,15){\vector(1,0){7.5}}
 \put(43,15){\vector(1,0){11.5}}
 \put(45,25){\vector(1,0){9.5}}
 \put(35,5){\vector(0,1){9.5}}
 \put(45,5){\vector(0,1){7.5}}
 \put(45,13){\vector(0,1){11.5}}
 \put(55,15){\vector(0,1){9.5}}
\put(32,3){$u$}
\put(39,1.5){1}
\put(32,9){1}
\put(46,7.5){0}
\put(37.5,16){0}
\put(49,16){1}
\put(42,19){1}
\put(56,19){0}
\put(49,26){0}
\put(85,25){\circle{1.0}}
\put(95,15){\circle{1.0}}
\put(95,25){\circle{1.0}}
\put(85,25){\vector(1,0){9.5}}
\put(95,15){\vector(0,1){9.5}}
\put(96,26){$v$}
\put(89,26){0}
\put(96,19){0}
\put(100,14){\line(1,0){9}}
\put(100,16){\line(1,0){9}}
\put(106,11){\line(1,1){4}}
\put(106,19){\line(1,-1){4}}
\put(115,5){\circle{1.0}}
\put(125,5){\circle{1.0}}
\put(115,15){\circle{1.0}}
\put(125,13){\circle{1.0}}
\put(123,15){\circle{1.0}}
\put(135,15){\circle{1.0}}
\put(125,25){\circle{1.0}}
\put(135,25){\circle{1.0}}
\put(115,5){\circle{2.0}}
\put(125,13){\circle{2.0}}
\put(123,15){\circle{2.0}}
\put(135,25){\circle{2.0}}
\put(115,5){\vector(1,0){9.5}}
\put(115,15){\vector(1,0){7.5}}
\put(123,15){\vector(1,0){11.5}}
\put(125,25){\vector(1,0){9.5}}
\put(115,5){\vector(0,1){9.5}}
\put(125,5){\vector(0,1){7.5}}
\put(125,13){\vector(0,1){11.5}}
\put(135,15){\vector(0,1){9.5}}
\put(136,26){$v$}
\put(119,1.5){1}
\put(112,9){1}
\put(126,7.5){0}
\put(117.5,16){0}
\put(129,16){1}
\put(122,19){1}
\put(136,19){0}
\put(129,26){0}
  \end{picture}
 \end{center}

Again, one shows that the label w.r.t. $i,j$ of each of the eight involved
edges is determined uniquely, just as indicated in the above picture (where the
bigger circles indicate critical vertices).

The final axiom concerns non-neighboring colors.

 \begin{itemize}
 \item[(A5)]
Let $|i-j|\ge 2$. Then for any $F\in\{F_i,F_i^{-1}\}$ and $F'\in
\{F_j,F_j^{-1}\}$, the operators $F,F'$ commute at each vertex where both
act.
  \end{itemize}

This is equivalent to saying that each component of the two-colored subgraph
$(V(K),E_i\sqcup E_j)$ is the Cartesian product of an $i$-path $P$ and a
$j$-path $P'$, or that each subcrystal of $K$ with non-neighboring colors $i,j$
is an $A_1\times A_1$-{\em crystal}.

 \medskip
One shows that any $A_n$-crystal $K$ is finite and has exactly one
zero-indegree vertex $s_K$ and one zero-outdegree vertex $t_K$, called the {\em
source} and {\em sink} of $K$, respectively. Furthermore, the $A_n$-crystals
$K$ admit a nice parameterization: the lengths $h_1(s_K),\ldots,h_n(s_K)$ of
monochromatic paths from the source determine $K$, and for each tuple
$c=(c_1,\ldots,c_n)$ of nonnegative integers, there exists a (unique)
$A_n$-crystal $K$ such that $c_i=h_i(s_K)$ for $i=1,\ldots,n$. (See~\cite{Stem}
and~\cite{cross}.) We call $c$ the {\em parameter} of $K$ and denote $K$ by
$K(c)$.

\subsection{The crossing model for $A_n$-crystals} \label{ssec:cross}

Following~\cite{cross}, the {\em crossing model} $\Mscr_n(c)$ generating the
$A_n$-crystal $K=K(c)$ with a parameter $c=(c_1,\ldots,c_n)\in \Zset_+^n$
consists of three ingredients:

(i) a directed graph $G_n=G=(V(G),E(G))$ depending on $n$, called the {\em
supporting graph} of the model;

(ii) a set $\Fscr=\Fscr(c)$ of {\em feasible} functions on $V(G)$;

(iii) a set $\Escr=\Escr(c)$ of transformations $f\mapsto f'$ of feasible
functions, called {\em moves}.

 \smallskip
To explain the construction of the supporting graph $G$, we first introduce
another directed graph $\Gscr=\Gscr_n$ that we call the {\em proto-graph} of
$G$. Its node set consists of elements $V_i(j)$ for all $i,j\in\{1,\ldots,n\}$
such that $j\le i$. (We use the term ``node'' for vertices in the crossing
model, to distinguish between these and vertices of crystals.) Its edges are
all possible pairs of the form $(V_i(j),V_{i-1}(j))$ ({\em ascending} edges) or
$(V_i(j),V_{i+1}(j+1))$ ({\em descending} edges). We say that the nodes
$V_i(1),\ldots,V_i(i)$ form $i$-th {\em level} of $\Gscr$ and order them as
indicated (by increasing $j$). We visualize $\Gscr$ by drawing it on the plane
so that the nodes of the same level lie in a horizontal line, the ascending
edges point North-East, and the descending edges point South-East. See the
picture where $n=4$.
  \begin{center}
  \unitlength=1mm
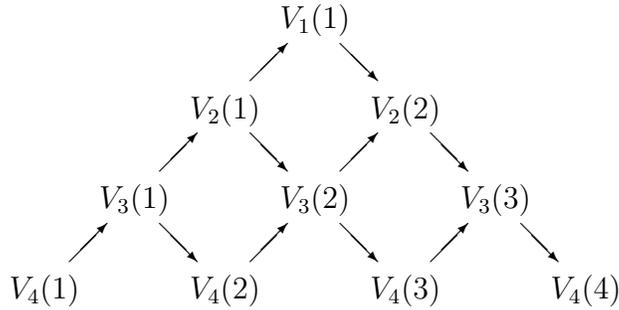
  \begin{picture}(80,40)
   \put(0,0){$V_4(1)$}
   \put(24,0){$V_4(2)$}
   \put(48,0){$V_4(3)$}
   \put(72,0){$V_4(4)$}
   \put(12,12){$V_3(1)$}
   \put(36,12){$V_3(2)$}
   \put(60,12){$V_3(3)$}
   \put(24,24){$V_2(1)$}
   \put(48,24){$V_2(2)$}
   \put(36,36){$V_1(1)$}
  \put(8,5){\vector(1,1){5}}
  \put(32,5){\vector(1,1){5}}
  \put(56,5){\vector(1,1){5}}
  \put(20,17){\vector(1,1){5}}
  \put(44,17){\vector(1,1){5}}
  \put(32,29){\vector(1,1){5}}
  \put(20,10){\vector(1,-1){5}}
  \put(44,10){\vector(1,-1){5}}
  \put(68,10){\vector(1,-1){5}}
  \put(32,22){\vector(1,-1){5}}
  \put(56,22){\vector(1,-1){5}}
  \put(44,34){\vector(1,-1){5}}
  \end{picture}
 \end{center}

The supporting graph $G$ is produced by replicating elements of $\Gscr$ as
follows. Each node $V_i(j)$ generates $n-i+1$ nodes of $G$, denoted as
$v_i^k(j)$ for $k=i-j+1,\ldots n-j+1$, which are ordered by increasing $k$ (and
accordingly follow from left to right in the visualization). We identify
$V_i(j)$ with the set of these nodes and call it a {\em multinode} of $G$. Each
edge of $\Gscr$ generates a set of edges of $G$ (a {\em multi-edge}) connecting
elements with equal upper indices. More precisely, $(V_i(j),V_{i-1}(j))$
produces $n-i+1$ ascending edges $(v_i^k(j),v_{i-1}^k(j))$ for $k=i-j+1,\ldots,
n-j+1$, and $(V_i(j),V_{i+1}(j+1))$ produces $n-i$ descending edges
$(v_i^k(j),v_{i+1}^k(j+1))$ for $k=i-j+1,\ldots,n-j$.

The resulting $G$ is the disjoint union of $n$ directed graphs
$G^1,\ldots,G^n$, where each $G^k$ contains all vertices of the form
$v_i^k(j)$. Also $G^k$ is isomorphic to the Cartesian product of two paths,
with the lengths $k-1$ and $n-k$. For example, for $n=4$, the graph $G$ is
viewed as
 \begin{center}
  \unitlength=1mm
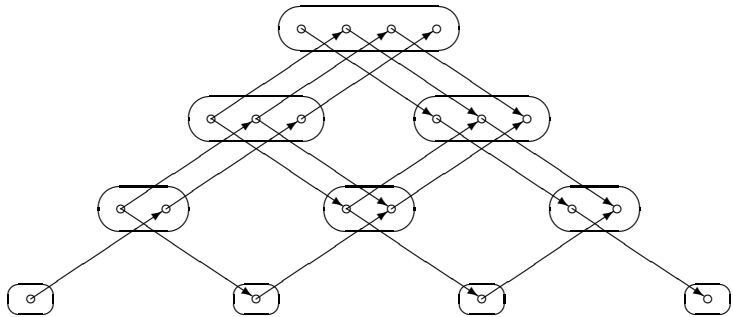
  \begin{picture}(90,39)
\put(0,0){\begin{picture}(54,36)%
\put(0,0){\circle{1.0}}
\put(18,12){\circle{1.0}}
\put(36,24){\circle{1.0}}
\put(54,36){\circle{1.0}}
\put(0,0){\vector(3,2){17.5}}
\put(18,12){\vector(3,2){17.5}}
\put(36,24){\vector(3,2){17.5}}
  \end{picture}}
\put(12,0){\begin{picture}(54,36)%
\put(0,12){\circle{1.0}}
\put(18,0){\circle{1.0}}
\put(18,24){\circle{1.0}}
\put(36,12){\circle{1.0}}
\put(36,36){\circle{1.0}}
\put(54,24){\circle{1.0}}
\put(0,12){\vector(3,2){17.5}}
\put(0,12){\vector(3,-2){17.5}}
\put(18,0){\vector(3,2){17.5}}
\put(18,24){\vector(3,2){17.5}}
\put(18,24){\vector(3,-2){17.5}}
\put(36,12){\vector(3,2){17.5}}
\put(36,36){\vector(3,-2){17.5}}
  \end{picture}}
\put(24,0){\begin{picture}(54,36)%
\put(0,24){\circle{1.0}}
\put(18,12){\circle{1.0}}
\put(18,36){\circle{1.0}}
\put(36,0){\circle{1.0}}
\put(36,24){\circle{1.0}}
\put(54,12){\circle{1.0}}
\put(0,24){\vector(3,2){17.5}}
\put(0,24){\vector(3,-2){17.5}}
\put(18,12){\vector(3,2){17.5}}
\put(18,12){\vector(3,-2){17.5}}
\put(18,36){\vector(3,-2){17.5}}
\put(36,0){\vector(3,2){17.5}}
\put(36,24){\vector(3,-2){17.5}}
  \end{picture}}
\put(36,0){\begin{picture}(54,36)%
\put(0,36){\circle{1.0}}
\put(18,24){\circle{1.0}}
\put(36,12){\circle{1.0}}
\put(54,0){\circle{1.0}}
\put(0,36){\vector(3,-2){17.5}}
\put(18,24){\vector(3,-2){17.5}}
\put(36,12){\vector(3,-2){17.5}}
  \end{picture}}
\put(45,36){\oval(24,6)}
\put(30,24){\oval(18,6)}
\put(60,24){\oval(18,6)}
\put(15,12){\oval(12,6)}
\put(45,12){\oval(12,6)}
\put(75,12){\oval(12,6)}
\put(0,0){\oval(6,4)}
\put(30,0){\oval(6,4)}
\put(60,0){\oval(6,4)}
\put(90,0){\oval(6,4)}
\end{picture}
 \end{center}
 \noindent
(where the multinodes are surrounded by ovals) and its components
$G^1,G^2,G^3,G^4$ are viewed as
 \begin{center}
  \unitlength=1mm
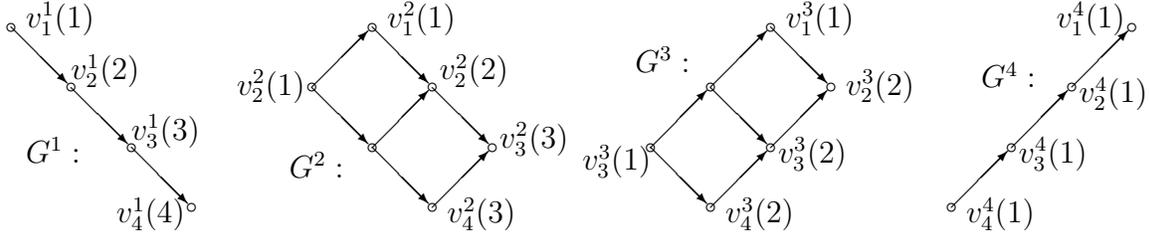
  \begin{picture}(150,27)
\put(0,24){\circle{1.0}}
\put(8,16){\circle{1.0}}
\put(16,8){\circle{1.0}}
\put(24,0){\circle{1.0}}
\put(0,24){\vector(1,-1){7.5}}
\put(8,16){\vector(1,-1){7.5}}
\put(16,8){\vector(1,-1){7.5}}
\put(2,6){$G^1:$}
\put(2,24){$v_1^1(1)$}
 \put(8,17){$v_2^1(2)$}
 \put(16,9){$v_3^1(3)$}
\put(14,-2){$v_4^1(4)$}
\put(40,16){\circle{1.0}}
\put(48,8){\circle{1.0}}
\put(48,24){\circle{1.0}}
\put(56,0){\circle{1.0}}
\put(56,16){\circle{1.0}}
\put(64,8){\circle{1.0}}
\put(40,16){\vector(1,1){7.5}}
\put(40,16){\vector(1,-1){7.5}}
\put(48,8){\vector(1,1){7.5}}
\put(48,8){\vector(1,-1){7.5}}
\put(48,24){\vector(1,-1){7.5}}
\put(56,0){\vector(1,1){7.5}}
\put(56,16){\vector(1,-1){7.5}}
\put(37,4){$G^2:$}
 \put(30,15){$v_2^2(1)$}
 \put(65,8){$v_3^2(3)$}
 \put(50,24){$v_1^2(1)$}
 \put(57,17){$v_2^2(2)$}
 \put(58,-2){$v_4^2(3)$}
\put(85,8){\circle{1.0}}
\put(93,0){\circle{1.0}}
\put(93,16){\circle{1.0}}
\put(101,8){\circle{1.0}}
\put(101,24){\circle{1.0}}
\put(109,16){\circle{1.0}}
\put(85,8){\vector(1,1){7.5}}
\put(85,8){\vector(1,-1){7.5}}
\put(93,0){\vector(1,1){7.5}}
\put(93,16){\vector(1,1){7.5}}
\put(93,16){\vector(1,-1){7.5}}
\put(101,8){\vector(1,1){7.5}}
\put(101,24){\vector(1,-1){7.5}}
\put(83,18){$G^3:$}
 \put(76,5){$v_3^3(1)$}
\put(95,-2){$v_4^3(2)$}
 \put(102,6){$v_3^3(2)$}
 \put(103,24){$v_1^3(1)$}
\put(111,15){$v_2^3(2)$}

\put(125,0){\circle{1.0}}
\put(133,8){\circle{1.0}}
\put(141,16){\circle{1.0}}
\put(149,24){\circle{1.0}}
\put(125,0){\vector(1,1){7.5}}
\put(133,8){\vector(1,1){7.5}}
\put(141,16){\vector(1,1){7.5}}
\put(129,16){$G^4:$}
 \put(127,-2){$v_4^4(1)$}
 \put(134,6){$v_3^4(1)$}
 \put(142,14){$v_2^4(1)$}
 \put(139,24){$v_1^4(1)$}
\end{picture}
 \end{center}

So each node $v=v_i^k(j)$ of $G$ has at most four incident edges, namely,
$(v_{i-1}^k(j-1),v)$, $(v_{i+1}^k(j),v)$, $(v,v_{i-1}^k(j))$,
$(v,v_{i+1}^k(j+1))$; we refer to them, when exist, as the NW-, SW-, NE-,
and SE-{\em edges}, and denote by $\eNW(v),\eSW(v),\eNE(v),\eSE(v)$,
respectively.

  \medskip
By a {\em feasible} function in the model (with a given $c$) we mean a function
$f:V(G)\to \Zset_+$ satisfying the following three conditions, where for an
edge $e=(u,v)$, $\;\partial f(e)$ denotes the increment $f(u)-f(v)$ of $f$ on
$e$, and $e$ is called {\em tight} for $f$, or $f$-{\em tight}, if $\partial
f(e)=0$:
  \begin{numitem1}
 \begin{itemize}
 \item[(i)] $f$ is {\em monotone} on the edges, in the sense that
$\partial f(e)\ge 0$ for all $e\in E(G)$;
 \item[(ii)] $0\le f(v)\le c_k$ for each $v\in V(G^k)$, $k=1,\ldots,n$;
 \item[(iii)] each multinode $V_i(j)$ contains a node $v$ with the following
property: the edge $\eSE(u)$ is tight for each node $u\in V_i(j)$
preceding $v$, and $\eSW(u')$ is tight for each node $u'\in V_i(j)$
succeeding $v$.
  \end{itemize}
 \label{eq:feas}
  \end{numitem1}

  \noindent
The {\em first} node $v=v_i^k(j)$ (i.e. with $k$ minimum) satisfying the
property in (iii) is called the {\em switch-node} of the multinode $V_i(j)$.
These nodes play an important role in our transformations of feasible functions
in the model.

To describe the rule of transforming $f\in\Fscr(c)$, we first extend each $G^k$
by adding extra nodes and edges (following~\cite{cross} and aiming to slightly
simplify the description). In the extended directed graph $\bar G^k$, the node
set consists of elements $v^k_i(j)$ for all $i=0,\ldots,n+1$ and $j=0,\ldots,n$
such that $j\le i$. The edge set of $\bar G^k$ consists of all possible pairs
of the form $(v^k_i(j),v^k_{i-1}(j))$ or $(v^k_i(j),v^k_{i+1}(j+1))$. Then all
$\bar G^k$ are isomorphic. The disjoint union of these $\bar G^k$ gives the
{\em extended supporting graph} $\bar G$. The creation of $\bar G^2$ from $G^2$
for $n=4$ is illustrated in the picture:

  \begin{center}
  \unitlength=1mm
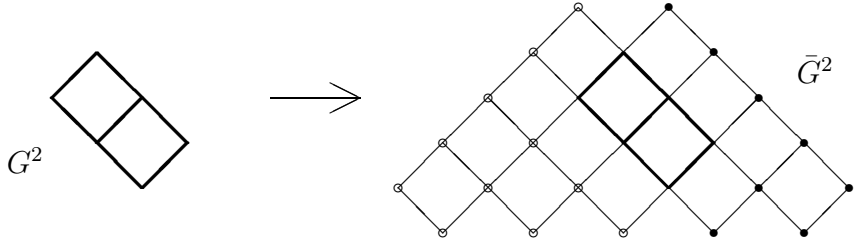
  \begin{picture}(120,30)
 \put(52,6){\line(1,1){24}}
 \put(58,0){\line(1,1){30}}
 \put(70,0){\line(1,1){24}}
 \put(82,0){\line(1,1){18}}
 \put(94,0){\line(1,1){12}}
 \put(106,0){\line(1,1){6}}
 \put(58,0){\line(-1,1){6}}
 \put(70,0){\line(-1,1){12}}
 \put(82,0){\line(-1,1){18}}
 \put(94,0){\line(-1,1){24}}
 \put(106,0){\line(-1,1){30}}
 \put(112,6){\line(-1,1){24}}
 \put(35,18){\line(1,0){12}}
 \put(47,18){\line(-2,1){4}}
 \put(47,18){\line(-2,-1){4}}
 \put(52,6){\circle{1.2}}
 \put(58,12){\circle{1.2}}
 \put(64,18){\circle{1.2}}
 \put(70,24){\circle{1.2}}
 \put(76,30){\circle{1.2}}
 \put(58,0){\circle{1.2}}
 \put(64,6){\circle{1.2}}
 \put(70,12){\circle{1.2}}
 \put(88,30){\circle*{1.2}}
 \put(70,0){\circle{1.2}}
 \put(76,6){\circle{1.2}}
 \put(94,24){\circle*{1.2}}
 \put(82,0){\circle{1.2}}
 \put(100,18){\circle*{1.2}}
 \put(94,0){\circle*{1.2}}
 \put(100,6){\circle*{1.2}}
 \put(106,12){\circle*{1.2}}
 \put(106,0){\circle*{1.2}}
 \put(112,6){\circle*{1.2}}
\thicklines{
 \put(6.1,18){\line(1,1){6}}
 \put(12.1,12){\line(1,1){6}}
 \put(18.1,6){\line(1,1){6}}
 \put(18.1,6){\line(-1,1){12}}
 \put(24.1,12){\line(-1,1){12}}
 \put(76.1,18){\line(1,1){6}}
 \put(82.1,12){\line(1,1){6}}
 \put(88.1,6){\line(1,1){6}}
 \put(88.1,6){\line(-1,1){12}}
 \put(94.1,12){\line(-1,1){12}}
 \put(5.9,18){\line(1,1){6}}
 \put(11.9,12){\line(1,1){6}}
 \put(17.9,6){\line(1,1){6}}
 \put(17.9,6){\line(-1,1){12}}
 \put(23.9,12){\line(-1,1){12}}
 \put(75.9,18){\line(1,1){6}}
 \put(81.9,12){\line(1,1){6}}
 \put(87.9,6){\line(1,1){6}}
 \put(87.9,6){\line(-1,1){12}}
 \put(93.9,12){\line(-1,1){12}}
 }
\put(0,8){$G^2$} \put(105,20){$\bar G^2$}
  \end{picture}
 \end{center}

Each feasible function on $V(G)$ is extended to the extra nodes $v=v^k_i(j)$ as
follows: $f(v):=c_k$ if there is a path from $v$ to a node of $G^k$, and
$f(v):=0$ otherwise (one may say that $v$ lies on the left of $G^k$ in the
former case, and on the right of $G^k$ in the latter case; in the above
picture, such nodes $v$ are marked by white and black circles, respectively).
In particular, each edge $e$ of $\bar G$ not incident with a node of $G$ is
tight, i.e. $\partial f(e)=0$ (extending $\partial f$ to the extra edges). For
a node $v=v_i^k(j)$ with $1\le j\le i\le n$ (and the given $f$), define the
value $\eps(v)=\eps_f(v)$ by
  \begin{equation}    \label{eq:epsv}
\eps(v):= \partial f(\eNW(v))-\partial f(\eSE(u))\quad
  (=\partial f(\eSW(v))-\partial f(\eNE(u)),
   \end{equation}
where $u:=v_i^k(j-1)$. For a multinode $V_i(j)$, define the numbers:
  \begin{eqnarray}
  \eps_i(j) &:=& \sum\left(\vphantom{\sum}\eps(v)\colon v\in V_i(j)\right); \label{eq:epsij} \\
    \eps_i(p,j) &:=& \eps_i(p)+\eps_i(p+1)+\ldots+\eps_i(j) \qquad
    \mbox{for}\;\;     1\le p \le j;    \nonumber \\
    \tilde\eps_i(j) &:=& \max\{0,\min\{\eps_i(p,j)\colon p=1,\ldots,j\}\}.
     \label{eq:barepsij}
  \end{eqnarray}
We call $\eps(v)$, $\eps_i(j)$ and $\tilde\eps_i(j)$ the {\em slack} at a node
$v$, the {\em total slack} at a multinode $V_i(j)$ and the {\em reduced slack}
at $V_i(j)$, respectively. (Note that $\eps,\tilde\eps$ are defined
in~\refeq{epsv},\refeq{epsij},\refeq{barepsij} in a somewhat different way than
in~\cite{cross}, which, however, does not affect the choice of active
multinodes and switch-nodes below.)

Now we are ready to define the transformations of $f$ (or the moves from $f$).
At most $n$ transformations $\phi_1,\ldots,\phi_n$ are possible. Each $\phi_i$
changes $f$ within level $i$ and is applicable when this level contains a
multinode $V_i(j')$ with $\tilde\eps_i(j')>0$. In this case we take the
multinode $V_i(j)$ such that
  \begin{equation} \label{eq:activej}
  \tilde\eps_i(j)>0\quad \mbox{and}\quad \tilde\eps_i(q)=0\;\;
                \mbox{for $q>j$},
  \end{equation}
referring to it as the {\em active} multinode for the given $f$ and $i$,
and increase $f$ by 1 at the
{\em switch-node} in $V_i(j)$, preserving $f$ on the other nodes of $G$. It is
shown~\cite{cross} that the resulting function $\phi_i(f)$ is again
feasible.

As a result, the model generates $n$-colored directed graph $\Kscr(c)=(\Fscr,
\Escr_1\sqcup\ldots\sqcup\Escr_n)$, where each color class $\Escr_i$ is formed
by the edges $(f,\phi_i(f))$ for all feasible functions $f$ to which the
operator $\phi_i$ is applicable. This graph is just an $A_n$-crystal.

  \begin{theorem} {\rm \cite[Th.~5.1]{cross}}  \label{tm:cross}
For each $n$ and $c\in\Zset_+^n$, the $n$-colored graph $\Kscr(c)$ is exactly
the $A_n$-crystal $K(c)$.
  \end{theorem}

 \subsection{Principal lattice and $(n-1)$-colored subcrystals of an $A_n$-crystal}
          \label{ssec:pr_lat}

Based on the crossing model, \cite{cross} reveals some important ingredients
and relations for an $A_n$-crystal $K=K(\bfc)$. One of them is the so-called
principal lattice, which is defined as follows.

Let $\bfa\in\Zset_+^n$ and $\bfa\le \bfc$. One easily checks that the function
on the vertices of the supporting graph $G$ that takes the constant value $a_k$
within each subgraph $G^k$ of $G$, $k=1,\ldots,n$, is feasible. We denote this
function and the vertex of $K$ corresponding to it by $f[\bfa]$ and
$\prv[\bfa]$, respectively, and call them {\em principal} ones. So the set of
principal vertices is bijective to the integer box
$\Bscr(\bfc):=\{a\in\Zset^n\colon \bfzero\le\bfa\le\bfc\}$; this set is called
the {\em principal lattice} of $K$ and denoted by $\Pi=\Pi(\bfc)$. When it is
not confusing, the term ``principal lattice'' may also be applied to
$\Bscr(\bfc)$.

The following properties of the principal lattice will be essentially used
later.

 \begin{prop} {\rm\cite[Expression~(6.4)]{cross}} \label{pr:fund_string}
Let $\bfa\in \Bscr(\bfc)$, $k\in\{1,\ldots,n\}$, and $\bfa':=\bfa+1_k$ (where
$1_k$ is $i$-th unit base vector in $\Rset^n$). The principal vertex
$\prv[\bfa']$ is obtained from $\prv[\bfa]$ by applying the operator string
   \begin{equation} \label{eq:string}
S_{n,k}:=w_{n,k,n-k+1}\cdots w_{n,k,2}w_{n,k,1},
   \end{equation}
where for $j=1,\ldots,n-k+1$, the substring $w_{n,k,j}$ is defined as
  $$
  w_{n,k,j}:=F_jF_{j+1}\cdots F_{j+k-1}.
  $$
When acting on $\Pi$, any two (applicable) strings $S_{n,k},S_{n,k'}$ commute.
In particular, any principal vertex $\prv[\bfa]$ is expressed via the source
$s_K=\prv[\bfzero]$ as
   \begin{equation} \label{eq:prin_str}
   \prv[\bfa]=S_{n,n}^{a_n}S_{n,n-1}^{a_{n-1}}\cdots S_{n,1}^{a_1}(s_K).
   \end{equation}
   \end{prop}

  \begin{prop} {\rm\cite[Prop.~6.1]{cross}} \label{pr:int_prlat}
For $c',c''\in Z_+^n$ with $c'\le c''\le c$, let $K(c'\colon\!c'')$ be the
subgraph of $K(c)$ formed by the vertices and edges contained in (directed)
paths from $\prv[c']$ to $\prv[c'']$ (the \emph{interval} of $K(c)$ from
$\prv[c']$ to $\prv[c'']$). Then $K(c'\colon\!c'')$ is isomorphic to the
$A_n$-crystal $K(c''-c')$, and the principal lattice of $K'$ consists of the
principal vertices $\prv[a]$ of $K(c)$ with $c'\le a\le c''$.
  \end{prop}

Let $\Kmn(c)$ denote the set of subcrystals with colors $1,\ldots,n-1$, and
$\Kmone$ the set of subcrystals with colors $2,\ldots,n$ in $K$ (recall that a
subcrystal is assumed to be connected and maximal for the corresponding subset
of colors).

  \begin{prop} {\rm\cite[Prop.~7.1]{cross}} \label{pr:prlat-subcryst}
Each subcrystal in $\Kmn$ (in $\Kmone$) contains precisely one principal
vertex. This gives a bijection between $\Kmn$ and $\Pi$ (resp., between
$\Kmone$ and $\Pi$).
  \end{prop}

We refer to the members of $\Kmn$ and $\Kmone$ as {\em upper} and {\em lower}
($(n-1)$-colored) subcrystals of $K$, respectively. For $\bfa\in \Bscr(\bfc)$,
the upper subcrystal containing the vertex $\prv[\bfa]$ is denoted by
$\Kup[\bfa]$. This subcrystal has its own principal lattice of dimension $n-1$,
which is denoted by $\Piup[\bfa]$. We say that the coordinate tuple $a$ is the
{\em locus} of $\Kup[\bfa]$ (and of $\Piup[\bfa]$) in $\Pi$. Analogously, for
$\bfb\in \Bscr(\bfc)$, the lower subcrystal containing $\prv[\bfb]$ is denoted
by $\Klow[\bfb]$, and its principal lattice by $\Pilow[\bfb]$; we say that
$\bfb$ is the locus of $\Klow[b]$ (and of $\Pilow[b]$) in $\Pi$. It turns out
that the parameters of upper and lower subcrystals can be expressed explicitly,
as follows.

  \begin{prop} {\rm\cite[Props.~7.2,7.3]{cross}} \label{pr:par-subcryst}
For $\bfa\in \Bscr(\bfc)$, the upper subcrystal $\Kup[\bfa]$ is isomorphic to
the $A_{n-1}$-crystal $K(\parup)$, where $\parup$ is the tuple
$(\parup_1,\ldots,\parup_{n-1})$ defined by
  \begin{equation} \label{eq:par_up}
  \parup_i:=c_i-a_i+a_{i+1}, \qquad i=1,\ldots,n-1.
  \end{equation}
The principal vertex $\prv[\bfa]$ is contained in the upper lattice
$\Piup[\bfa]$ and its coordinate $\heartup=(\heartup_1,\ldots,\heartup_{n-1})$
in $\Piup[\bfa]$ satisfies
  \begin{equation}  \label{eq:heart_up}
  \heartup_i=a_{i+1}, \qquad i=1,\ldots,n-1.
  \end{equation}

Symmetrically, for $\bfb\in \Bscr(\bfc)$, the lower subcrystal $\Klow[\bfb]$ is
isomorphic to the $A_{n-1}$-crystal $K(\parlow)$ with colors $2,\ldots,n$,
where $\parlow$ is defined by
  \begin{equation} \label{eq:par_low}
  \parlow_i:=c_i-b_i+b_{i-1}, \qquad i=2,\ldots,n.
  \end{equation}
The principal vertex $\prv[\bfb]$ is contained in the lower lattice $\Pilow[b]$
and its coordinate $\heartlow=(\heartlow_2,\ldots,\heartlow_n)$ in $\Pilow[b]$
satisfies
  \begin{equation}  \label{eq:heart_low}
  \heartlow_i=b_{i-1}, \qquad i=2,\ldots,n.
  \end{equation}
  \end{prop}

We call $\prv[\bfa]$ the {\em heart} of $\Kup[\bfa]$ w.r.t. $K$, and similarly
for lower subcrystals.

(As one more result, \cite{cross} also gives (in Remark~5) a piecewise linear
formula to compute, given an $(n-1)$-tuple $q$, the number of upper subcrystals
of $K(c)$ with the parameter equal to $q$, but we do not need this in what
follows.)
 \smallskip

 \noindent
{\bf Remark 1.} As is mentioned in the Introduction, the crossing model is
viewed as a refinement of the Gelfand-Tsetlin pattern (or \emph{GT-pattern})
model~\cite{GT-50}. More precisely, for $c\in\Zset_+^n$, form the partition
$\lambda=(\lambda_1\ge\lambda_2\ge\cdots\ge\lambda_n\ge\lambda_{n+1}=0)$ by
setting $\lambda_i:=c[1:n-i+1]$, where $c[p:q]$ denotes $c_p+c_{p+1}+\ldots
+c_q$. A GT-pattern for $\lambda$ is a triangular array $X=(x_{ij})_{1\le j\le
i\le n}$ of integers satisfying: (a) $x_{ij}\ge x_{i-1,j},x_{i+1,j+1}$, and (b)
$\lambda_i\ge x_{n,i}\ge \lambda_{i+1}$, for all possible $i,j$. As is
explained in~\cite{cross}, the set of feasible functions $f$ in the crossing
model $\Mscr_n(c)$ is bijective to the set of GT-patterns $X$ for $\lambda$;
such a correspondence is given by $x_{i,j}:=\hat f_i(j)+c[1:i-j]$, where $\hat
f_i(j)$ denotes the sum of values of $f$ over the multinode $V_i(j)$. However,
it is not clear how to visualize, and work with, principal vertices directly in
terms of GT-patterns, whereas such vertices are well visualized and convenient
to handle in the crossing model.


\section{Assembling an $A_n$-crystal} \label{sec:ass_A}

As mentioned in the Introduction, the structure of an $A_n$-crystal $K=K(\bfc)$
will be described in a recursive manner. The idea is as follows. We know that
$K$ contains $|\Pi|=(c_1+1)\times\ldots\times (c_n+1)$ upper subcrystals (with
colors $1,\ldots,n-1$) and $|\Pi|$ lower subcrystals (with colors
$2,\ldots,n$). Moreover, the parameters of these subcrystals are expressed
explicitly by~\refeq{par_up} and~\refeq{par_low}. Assume that the set
$\Kscr^{(-n)}$ of upper subcrystals and the set $\Kscr^{(-1)}$ of lower
subcrystals are available. Then, in order to assemble $K$, it suffices to point
out, in appropriate terms, the intersection $\Kup[\bfa]\cap \Klow[\bfb]$ for
all pairs $\bfa,\bfb\in \Bscr(\bfc)$ (the intersection may either be empty, or
consist of one or more $(n-2)$-colored subcrystals with colors $2,\ldots,n-1$
in $K$). We give an appropriate characterization in Theorems~\ref{tm:two_div}
and~\ref{tm:two_loci} below.

To state them, we need additional terminology and notation. Consider a
subcrystal $\Kup[\bfa]$, and let $\parup,\heartup$ be defined as
in~\refeq{par_up},\refeq{heart_up}. For $\bfp=(p_1,\ldots,p_{n-1})\in
\Bscr(\parup)$, the vertex in the upper lattice $\Piup[\bfa]$ having the
coordinate $\bfp$ is denoted by $\vup[\bfa,\bfp]$. We call the vector
$\Delta:=\bfp-\heartup$ the {\em deviation} of $\vup[\bfa,\bfp]$ from the heart
$\prv[\bfa]$ in $\Piup[\bfa]$, and will use the alternative notation
$\vup[a|\Delta]$ for this vertex. In particular,
$\prv[a]=\vup[\bfa,\heartup]=\vup[\bfa|\,0]$.

Similarly, for a lower subcrystal $\Klow[\bfb]$, let $\parlow,\heartlow$ be as
in~\refeq{par_low},\refeq{heart_low}. For $\bfq=(q_2,\ldots,q_n)\in
\Bscr(\parlow)$, the vertex with the coordinate $\bfq$ in $\Pilow[\bfb]$ is
denoted by $\vlow[\bfa,\bfq]$. Its deviation is $\nabla:=\bfq-\heartlow$, and
we may alternatively denote this vertex by $\vlow[b|\nabla]$.

We call an $(n-2)$-colored subcrystal with colors $2,\ldots,n-1$ in $K$ a {\em
middle subcrystal} and denote the set of these by $\Kscr^{(-1,-n)}$. Each
middle crystal $\Kmid$ is a {\em lower} subcrystal of some upper subcrystal
$K'=\Kup[a]$ of $K$. By Proposition~\ref{pr:prlat-subcryst} applied to $K'$,
~$\Kmid$ has a unique vertex $\vup[a|\Delta]$ in the lattice $\Piup[\bfa]$. So
each $\Kmid$ can be encoded by a pair $(a,\Delta)$ formed by a point $a\in
\Bscr(c)$ and a deviation $\Delta$ in $\Piup[a]$. At the same time, $\Kmid$ is
an \emph{upper} subcrystal of some lower subcrystal $\Klow[b]$ of $K$ and has a
unique vertex $\vlow[b|\nabla]$ in $\Pilow[b]$. Therefore, the members of
$\Kscr^{(-1,-n)}$ determine a bijection
  $$
  \zeta:(a,\Delta)\mapsto (b,\nabla)
  $$
between all pairs $(a,\Delta)$ concerning upper subcrystals and all pairs
$(b,\nabla)$ concerning lower subcrystals.

The map $\zeta$ is expressed explicitly in the following two theorems. Here for
a tuple $\rho=(\rho_i\colon i\in I)$ of reals, we denote by $\rho^+$ ($\rho^-$)
the tuple formed by $\rho_i^+:=\max\{0,\rho_i\}$ (resp.
$\rho_i^-:=\min\{0,\rho_i\}$), $i\in I$.

 \begin{theorem} \label{tm:two_div} {\rm(on two deviations).}
Let $a\in \Bscr(c)$ and let $\Delta=(\Delta_1,\ldots,\Delta_{n-1})$ be a
deviation in $\Piup[a]$ (from the heart of $\Kup[a]$). Let
$(b,\nabla)=\zeta(a,\Delta)$. Then
  \begin{equation} \label{eq:nablai}
  \nabla_i=-\Delta_{i-1}, \qquad i=2,\ldots,n.
 \end{equation}
 \end{theorem}

  \begin{theorem} \label{tm:two_loci} {\rm(on two loci).}
Let $a,b,\Delta,\nabla$ be as in the previous theorem. Then
  \begin{equation} \label{eq:bi}
  b_i=a_i+\Delta^+_i+\Delta^-_{i-1}, \qquad i=1,\ldots,n,
  \end{equation}
letting $\Delta_0=\Delta_n:=0$.
 \end{theorem}

Proofs of these theorems will be given in the next section.

Based on Theorems~\ref{tm:two_div} and~\ref{tm:two_loci}, the crystal $K(c)$ is
assembled as follows. By recursion we assume that all upper and lower
subcrystals are already constructed. We also assume that for each upper
subcrystal $\Kup[a]\simeq K(\parup)$, its principal lattice is distinguished by
use of the corresponding injective map $\sigma:\Bscr(\parup)\to V(K(\parup))$,
and similarly for the lower subcrystals. We delete the edges with color 1 in
each $\Kup[a]$ and extract the components of the resulting graphs, forming the
set $\Kscr^{(-1,-n)}$ (arranged as a list) of all middle subcrystals of $K(c)$.
Each $\Kmid \in\Kscr^{(-1,-n)}$ is encoded by a corresponding pair
$(a,\Delta)$, where $\bfa\in\Bscr(\bfc)$ and the deviation $\Delta$ in
$\Piup[a]$ is determined by use of $\sigma$ as above. Acting similarly for the
lower subcrystals $\Klow[b]$ (by deleting the edges  with color $n$), we obtain
the same set of middle subcrystals (arranged as another list), each of which
being encoded by a corresponding pair $(b,\nabla)$, where $b\in\Bscr(c)$ and
$\nabla$ is a deviation in $\Pilow(b)$. Relations~\refeq{bi} and \refeq{nablai}
indicate how to identify each member of the first list with its counterpart in
the second one. Now restoring the deleted edges with colors 1 and $n$, we
obtain the desired crystal $K(c)$. The corresponding map $\Bscr(c)\to V(K(c))$
is constructed easily (e.g., by use of operator strings as in
Proposition~\ref{pr:fund_string}).

  \smallskip
We conclude this section with several remarks.
 \smallskip

 \noindent
{\bf Remark 2.} For each $a\in \Bscr(c)$ and each vertex $v=\vup[a,p]$ in the
upper lattice $\Piup[a]$, one can express the parameter $\parmid= (\parmid_2,
\ldots,\parmid_{n-1})$ of the middle subcrystal $\Kmid$ containing $v$, as well
as the coordinate $\heartmid=(\heartmid_2,\ldots, \heartmid_{n-1})$ of its
heart w.r.t. $\Kup[a]$ in the principal lattice of $\Kmid$. Indeed, since
$\Kmid$ is a lower subcrystal of $\Kup[a]$, one can apply relations as
in~\refeq{par_low},\refeq{heart_low}. Denoting the parameter of $\Kup[a]$ by
$\parup$ and the coordinate of its heart in $\Piup[a]$ by $\heartup$, letting
$\Delta:=p-\heartup$, and using~\refeq{par_up},\refeq{heart_up}, we have for
$i=2,\ldots,n-1$:
  \begin{eqnarray}
  &&\parmid_i =\parup_i-p_i+p_{i-1}=(c_i-a_i+a_{i+1})-(a_{i+1}+\Delta_i)
     +(a_i+\Delta_{i-1})\qquad \nonumber  \\
   & & \qquad\qquad\qquad\qquad\qquad\qquad\qquad\qquad\qquad\qquad
            =c_i-\Delta_i+\Delta_{i-1}; \label{eq:par_mid} \\
  &&\heartmid_i = p_{i-1}=\heartup_{i-1}+\Delta_{i-1}=a_i+\Delta_{i-1}.
  \nonumber
  \end{eqnarray}

Symmetrically, if $\Kmid$ is contained in $\Klow[b]$ and has deviation $\nabla$
in $\Pilow[b]$, then for $i=2,\ldots,n-1$:
   \begin{eqnarray}
   &&\parmid_i = c_i-\nabla_i+\nabla_{i+1}; \label{eq:par_mid2} \\
   && \bmid_i=b_i+\nabla_{i+1}, \nonumber
   \end{eqnarray}
where $\bmid$ is the coordinate of the heart of $\Kmid$ w.r.t. $\Klow[b]$ in
the principal lattice of $\Kmid$ (note that $\bmid$ may differ from
$\heartmid$). We will use~\refeq{par_mid} and~\refeq{par_mid2} in
Section~\ref{sec:proof}.

 \medskip
 \noindent
{\bf Remark 3.} A straightforward implementation of the above recursive method
of constructing $K=K(c)$ takes $O(2^{q(n)}N)$ time and space, where $q(n)$ is a
polynomial in $n$ and $N$ is the number of vertices of $K$. Here the factor
$2^{q(n)}$ appears because the total number of vertices in the upper and lower
subcrystals is $2N$ (implying that there appear $4N$ vertices in total on the
previous step of the recursion, and so on). Therefore, such an implementation
has polynomial complexity of the size of the output for each fixed $n$, but not
in general. However, many intermediate subcrystals arising during the recursive
process are repeated, and we can use this fact to improve the implementation.
More precisely, the colors occurring in each intermediate subcrystal in the
process form an interval of the ordered set $(1,\ldots,n)$. We call a
subcrystal of this sort a {\em color-interval subcrystal}, or a {\em
CI-subcrystal}, of $K$. In fact, {\em every} CI-subcrystal of $K$ appears in
the process. Since the number of intervals is $\frac{n(n+1)}{2}$ and the
CI-subcrystals concerning equal intervals are pairwise disjoint, the total
number of vertices of all CI-subcrystals of $K$ is $O(n^2 N)$. It is not
difficult to implement the recursive process so that each CI-subcrystal $K'$ be
explicitly constructed only once. As a result, we obtain the following
  \begin{prop} \label{pr:efficient}
Let $c\in\Zset_+^n$. The $A_n$-crystal $K(c)$ and all its CI-subcrystals can be
constructed in $O(q'(n)|V(K(c))|)$ time and space, where $q'(n)$ is a
polynomial in $n$.
  \end{prop}

 \noindent
{\bf Remark 4.} Relation~\refeq{bi} shows that the intersection of $\Kup[a]$
and $\Klow[b]$ may consist of many middle subcrystals. Indeed, if $\Delta_i>0$
and $\Delta_{i-1}<0$ for some $i$, then $b$ does not change by simultaneously
decreasing $\Delta_i$ by 1 and increasing $\Delta_{i-1}$ by 1. The number of
common middle subcrystals of $\Kup[a]$ and $\Klow[b]$ for arbitrary $a,b\in
\Bscr(c)$ can be expressed by an explicit piecewise linear formula,
using~\refeq{bi} and the box constraints $-a_{i+1}\le\Delta_i\le c_i-a_i$,
$i=1,\ldots,n-1$, on the deviations $\Delta$ in $\Piup[a]$ (which follow
from~\refeq{par_up},\refeq{heart_up}).

 \section{Proofs of Theorems~\ref{tm:two_div} and~\ref{tm:two_loci}} \label{sec:proof}

Let $a,\Delta,b,\nabla$ be as in the hypotheses of Theorem~\ref{tm:two_div}.
First we show that Theorem~\ref{tm:two_div} follows from
Theorem~\ref{tm:two_loci}.

 \medskip
 \noindent
{\bf Proof of~\refeq{nablai}} \emph{in the assumption that~\refeq{bi} is
valid}. The middle subcrystal $\Kmid$ determined by $(a,\Delta)$ is the same as
the one determined by $(b,\nabla)$. The parameter $\parmid$ of $\Kmid$ is
expressed simultaneously by~\refeq{par_mid} and by~\refeq{par_mid2}. Then
$c_i-\Delta_i+\Delta_{i-1}=c_i-\nabla_i+\nabla_{i+1}$ for $i=2,\ldots,n-1$.
Therefore,
   \begin{equation}  \label{eq:Delta-nabla}
  \Delta_1+\nabla_2=\Delta_2+\nabla_3=\ldots=
         \Delta_{n-1}+\nabla_n=:\alpha.
   \end{equation}

In order to obtain~\refeq{nablai}, one has to show that $\alpha=0$. We argue as
follows. Renumber the colors $1,\ldots,n$ as $n,\ldots,1$, respectively; this
yields the crystal $\hat K= K(\hat c)$ symmetric to $K(c)$. Then $\Klow[b]$
turns into the upper subcrystal $\hat K^\uparrow[\hat b]$ of $\hat K$, where
$(\hat b_1,\ldots,\hat b_n)=(b_n,\ldots,b_1)$. Also the deviation $\nabla$ in
$\Pilow[b]$ turns into the deviation $\hat\nabla=(\hat\nabla_1,
\ldots,\hat\nabla_{n-1})=(\nabla_n,\ldots,\nabla_2)$ in the principal lattice
of $\hat K^\uparrow[\hat b]$. Applying relations as in~\refeq{bi} to $(\hat
b,\hat\nabla)$, we have
  \begin{equation} \label{eq:bar_nabla}
  \hat a_i=\hat b_i+\hat\nabla^+_i+\hat\nabla^-_{i-1}=
    b_{n-i+1}+\nabla^+_{n-i+1}+\nabla^-_{n-i+2}, \qquad i=1,\ldots,n,
   \end{equation}
where $\hat a_i:=a_{n-i+1}$ and $\hat\nabla^+_n:=\hat\nabla^-_0:=0$.
On the other hand,~\refeq{bi} for $(a,\Delta)$ gives
  \begin{equation} \label{eq:nab}
  \hat b_i=b_{n-i+1}=a_{n-i+1}+\Delta^+_{n-i+1}+\Delta^-_{n-i},
            \qquad i=1,\ldots,n.
  \end{equation}

Relations~\refeq{bar_nabla} and~\refeq{nab} imply
  $$
  a_{n-i+1}=b_{n-i+1}+\nabla^+_{n-i+1}+\nabla^-_{n-i+2}=
      (a_{n-i+1}+\Delta^+_{n-i+1}+\Delta^-_{n-i})
    +\nabla^+_{n-i+1}+\nabla^-_{n-i+2},
    $$
whence
   $$
 \Delta^+_{n-i+1}+\Delta^-_{n-i}
    +\nabla^+_{n-i+1}+\nabla^-_{n-i+2}=0, \qquad i=1,\ldots,n.
  $$
Adding up the latter equalities, we obtain
  $$
  (\Delta_1+\ldots+\Delta_{n-1})+(\nabla_2+\ldots+\nabla_n)=0.
  $$
This and~\refeq{Delta-nabla} imply $(n-1)\alpha=0$. Hence $\alpha=0$,
yielding~\refeq{nablai} and Theorem~\ref{tm:two_div}. \hfill \qed
 \bigskip

 \noindent
{\bf Proof of Theorem~\ref{tm:two_loci}.} It is more intricate and essentially
uses the crossing model.

For a feasible function $f\in\Fscr(c)$ and its corresponding vertex $v$ in
$K=K(c)$, we may denote $v$ as $v_f$, and $f$ as $f_v$. From the crossing model
it is seen  that
  \begin{numitem1}
if a vertex $v\in V(K)$ belongs to $\Kup[a]$ and to $\Klow[b]$, then the
tuples $a$ and $b$ are expressed via the values of $f=f_v$ in levels $n$
and 1 as follows:

$a_k=f(v_n^k(n-k+1))$\; and\; $b_k=f(v_1^k(1))$\; for $k=1,\ldots,n$.
  \label{eq:fn1}
  \end{numitem1}
Indeed, the principal vertex $\prv[a]$ is reachable from $v$ by applying
operators $F_i$ or $F_i^{-1}$ with $i\ne n$. The corresponding moves in the
crossing model do not change $f$ within level $n$. Similarly, $\prv[b]$ is
reachable from $v$ by applying operators $F_i$ or $F_i^{-1}$ with $i\ne 1$, and
the corresponding moves do not change $f$ within level 1. Also the first
(second) equality in~\refeq{fn1} is valid for the principal function
$f_{\prv[a]}$ (resp. $f_{\prv[b]}$).

Next we introduce special functions on the node set $V(G)$ of the supporting
graph $G=G_n$. Consider a component $G^k=(V^k,E^k)$ of $G$. It is a rectangular
grid of size $k\times(n-k+1)$ (rotated by 45$^\circ$ in the visualization of
$G$), and its vertex set is
  $$
  V^k=\{v_i^k(j)\colon j=1,\ldots,n-k+1,\;\; i=j,\ldots,j+k-1\}.
  $$
To represent it in a more convenient form, let us rename $v_i^k(j)$ as
$x^k_{i-j+1}(j)$, or $x_{i-j+1}(j)$ (as though rotating $G^k$ by $45^\circ$).
Then
  $$
  V^k=\{x_m(j)\colon  m=1,\ldots,k,\;\; j=1,\ldots,n-k+1\},
  $$
the SE-edges of $G^k$ become of the form $(x_m(j),x_m(j+1))$, and the NE-edges
become of the form $(x_m(j),x_{m-1}(j))$. We distinguish the following subsets
of $V^k$: \smallskip

(i) the {\em SW-side} $P^k:=\{x_k(1),\ldots,x_k(n-k+1)\}$;

(ii) the {\em right rectangle} $R^k:=\{x_m(j)\colon 1\le m\le k-1,\; 1\le j\le
n-k+1\}$;

(iii) the {\em left rectangle} $L^k:=\{x_m(j)\colon 1\le m\le k,\; 1\le j\le
n-k\}$.

\noindent Denote the characteristic functions (in $\Rset^{V^k}$) of $P^k, R^k,
L^k$ as $\pi^k$, $\rho^k,\lambda^k$, respectively.

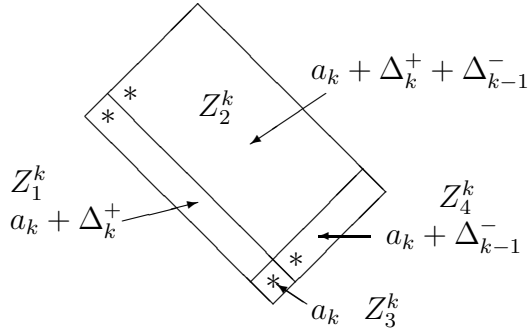
\begin{figure}[hbt]                  
  \begin{center}
  \unitlength=1mm
  \begin{picture}(70,40)(0,6)
  \put(10,30){\line(1,1){15}}
  \put(32,8){\line(1,1){15}}
  \put(35,5){\line(1,1){15}}
  \put(10,30){\line(1,-1){25}}
  \put(13,33){\line(1,-1){25}}
  \put(25,45){\line(1,-1){25}}
  \put(40,35){$a_k+\Delta^+_k+\Delta^-_{k-1}$}
  \put(12,29){$\ast$}
  \put(15,32){$\ast$}
  \put(34,7){$\ast$}
  \put(37,10){$\ast$}
  \put(0,20){$Z_1^k$}
  \put(0,15){$a_k+\Delta^+_k$}
  \put(40,3){$a_k$}
  \put(47,3){$Z^k_3$}
  \put(50,13){$a_k+\Delta^-_{k-1}$}
  \put(57,18){$Z^k_4$}
  \put(25,30){$Z^k_2$}
  \put(42,33){\vector(-3,-2){10}}
  \put(15,16){\vector(4,1){10}}
  \put(39.5,5){\vector(-3,2){4}}
  \put(48,14){\vector(-1,0){7}}
 \end{picture}
 \end{center}
  \caption{The partition of $V^k$.}
  \label{fig:partition}
  \end{figure}

Return to $a$ and $\Delta$ as above. We associate to $(a,\Delta)$ the functions
  \begin{equation} \label{eq:fkaD}
  f^k_{a,\Delta}:=a_k\pi^k+(a_k+\Delta^-_{k-1})\rho^k+\Delta^+_k\lambda^k
  \end{equation}
on $V^k$ for $k=1,\ldots,n$ (see Fig.~\ref{fig:partition}), and define
$f_{a,\Delta}$ to be the function on $V(G)$ whose restriction to each $V^k$ is
$f^k_{a,\Delta}$.

In view of~\refeq{fn1}, $f=f_{a,\Delta}$ takes the values in levels $n$ and 1
as required in~\refeq{bi} (with $k$ in place of $i$), namely,
$f(v^k_n(n-k+1))=a_k$ and $f(v^k_1(1))=a_k+\Delta^+_k+\Delta^-_{k-1}$ for
$k=1,\ldots,n$. Therefore, to obtain~\refeq{bi} it suffices to show the
following
  \begin{lemma} \label{lm:faDelta}
{\rm(i)} The function $f=f_{a,\Delta}$ is feasible. {\rm(ii)} The vertex $v_f$
belongs to $\Piup[a]$ and has the deviation $\Delta$ in it; in other words,
$f=f_{\prv[a\vert\Delta]}$.
  \end{lemma}

\noindent {\bf Proof}~ First we prove assertion~(i). Let $k\in\{1,\ldots,n\}$.
We partition $V^k$ into four subsets (rectangular {\em pieces}):
  $$
  Z^k_1:=P^k- \{x_k^k(n-k+1)\};\;\;
  Z^k_2:=L^k\setminus P^k;\;\;
  Z^k_3:=\{x_k^k(n-k+1)\};\;\;
  Z^k_4:=R^k\setminus L^k
  $$
\noindent (where $Z^k_2=Z^k_4=\emptyset$ when $k=1$, and
$Z^k_1=Z^k_2=\emptyset$ when $k=n$). By~\refeq{fkaD},
  \begin{numitem1}
$f$ takes a constant value within each piece $Z^k_q$, namely: $a_k+\Delta^+_k$
on $Z^k_1$; \\
$a_k+\Delta^+_k+\Delta^-_{k-1}$ on $Z^k_2$; $a_k$ on $Z^k_3$; and
$a_k+\Delta^-_{k-1}$ on $Z^k_4$
  \label{eq:valuesZ}
  \end{numitem1}
(as illustrated in Fig.~\ref{fig:partition}). Also each edge of $G^k$
connecting different pieces goes either from $Z^k_1$ to $Z^k_2\cup Z^k_3$ or
from $Z^k_2\cup Z^k_3$ to $Z^k_4$. This and~\refeq{valuesZ} imply that
$\partial f(e)\ge 0$ for each edge $e\in E^k$, whence $f$
satisfies~\refeq{feas}(i).

The deviation $\Delta$ is restricted as $-\heartup\le\Delta\le
\parup-\heartup$, where $\parup$ is the parameter of the subcrystal $\Kup[a]$
and $\heartup$ is the coordinate of its heart $\prv[a]$ in $\Piup[a]$.
Formulas~\refeq{par_up} and~\refeq{heart_up} for $\parup$ and $\heartup$ give
  \begin{equation} \label{eq:boundsDelta}
  -a_{k+1}\le\Delta_k\le c_k-a_k\quad \mbox{and}\quad
  -a_k\le\Delta_{k-1}\le c_{k-1}-a_{k-1}.
  \end{equation}
The inequalities $\Delta_k\le c_k-a_k$ and $a_k\le c_k$ imply
$a_k+\Delta^+_k\le c_k$. The inequalities $-a_k\le\Delta_{k-1}$ and $a_k\ge 0$
imply $a_k+\Delta^-_{k-1}\ge 0$. Then, in view of~\refeq{valuesZ}, we obtain
$0\le f(v)\le c_k$ for each node $v$ of $G^k$, yielding~\refeq{feas}(ii).

To verify the switch condition~\refeq{feas}(iii), consider a multinode
$V_i(j)$ with $i<n$. It consists of $n-i+1$ nodes $v^k_i(j)$, where
$i-j+1\le k\le n-j+1$.

Let $i\le n-2$. Suppose that $v=v_i^k(j)$ is a node whose SW-edge $e=(u,v)$
exists and is not $f$-tight. This is possible only if $u\in Z^k_1$ and $v\in
Z^k_2$. In this case $k$ is determined as $k=i-j+2$, i.e. $v$ is the second
node in $V_i(j)$. We observe that: (a) for the first node $v_i^{k-1}(j)$ of
$V_i(j)$, both ends of its SE-edge $e'$ belong to the piece $Z_1^{k-1}$; and
(b) for any node $v_i^{k'}(j)$ with $k'>k$ in $V_i(j)$, both ends of its
SW-edge $e''$ belong either to $Z_2^{k'}$ or to $Z_4^{k'}$. So such $e'$ and
$e''$ are $f$-tight. Therefore, the node $v$ satisfies the condition
in~\refeq{feas}(iii) for $V_i(j)$.

Now let $i=n-1$. Then $V_i(j)$ consists of two nodes $v=v^{n-j}_{n-1}(j)$ and
$v'=v^{n-j+1}_{n-1}(j)$. Put $k:=n-j$. Then the edge $e=\eSE(v)$ goes from
$Z^k_1$ to $Z^k_3=\{x^k_k(n-k+1)\}$, and the edge $e'=\eSW(v')$ goes from
$Z^{k+1}_3=\{x^{k+1}_{k+1}(n-k)\}$ to $Z^{k+1}_4$. By~\refeq{valuesZ}, we have
$\partial f(e)=(a_k+\Delta_k^+)-a_k =\Delta^+_k$ and $\partial
f(e')=a_{k+1}-(a_{k+1}+\Delta_k^-)=-\Delta^-_k$. Since at least one of
$\Delta^+_k,\Delta^-_k$ is zero, we conclude that at least one of $e,e'$ is
tight. So~\refeq{feas}(iii) is valid again.

 \medskip
Next we prove assertion~(ii) in the lemma. (The idea is roughly as follows. For
each $k$, compare the function $f^k_{a,\Delta}=:g$ with the function $h$ on
$V^k$ taking the constant value $a_k$. By~\refeq{fkaD}, $g=h+\Delta^+_k
\lambda^k+ \Delta^-_{k-1} \rho^k$. In other words, $g$ is obtained from
${f_{\prv[a]}}\rest{V^k}$ by adding $\Delta^+_k$ times the ``left rectangle
function'' $\lambda^k$, followed by subtracting $|\Delta^-_{k-1}|$ times the
``right rectangle function'' $\rho^k$. A crucial observation is that adding
$\lambda^k$ corresponds to applying the operator string $S_{n-1,k}$ (or
shifting by $k$-th unit base vector in the upper principal lattice $\Piup[a]$),
while subtracting $\rho^k$ corresponds to applying $S^{-1}_{n-1,k-1}$ (or
shifting by minus $(k-1)$-th unit base vector in $\Piup[a]$).  This is because
the substrings $w$ in $S_{n-1,k}$ correspond to the SW--NE paths in $L^k$, and
the substrings in $S_{n-1,k-1}$ to similar paths in $R^{k}\simeq L^{k-1}$.)

Now we give a more careful and formal description. We use induction on
  $$
  \eta(\Delta):=\Delta_1+\ldots+\Delta_{n-1}.
  $$

In view of~\refeq{boundsDelta}, $\eta(\Delta)\ge -a_2-\ldots-a_n$. Suppose this
turns into equality. Then $\Delta_k=-a_{k+1}\le 0$ for $k=1,\ldots,n-1$, and
$f=f_{a,\Delta}$ takes the following values within each $V^k$
(cf.~\refeq{valuesZ}): $f(v)=a_k$ if $v\in P^k$, and $f(v)=0$ if $v\in
V^k-P^k$. This $f$ is the minimal feasible function whose values in level $n$
correspond to $a$, i.e. $v_f$ is the source of $\Kup[a]$. Then $v_f$ is the
minimal vertex $\prv[a,\bfzero]$ in $\Piup[a]$, and its deviation in $\Piup[a]$
is just $\Delta$, as required. This gives the base of our induction.

Now consider an arbitrary $\Delta$ satisfying~\refeq{boundsDelta}. Let $k$ be
such that $\Delta_k<c_k-a_k$ (if any) and define $\Delta'_k:=\Delta_k+1$ and
$\Delta'_i:=\Delta_i$ for $i\ne k$. Then $\eta(\Delta)<\eta(\Delta')$. We
assume by induction that assertion~(ii) is valid for $f_{a,\Delta}$, and our
aim is to show validity of~(ii) for $f_{a,\Delta'}$.

In what follows $f$ stands for the former function $f_{a,\Delta}$.

Let $v'$ be the vertex with the deviation $\Delta'$ in $\Piup[a]$. Both $v_f$
and $v'$ are principal vertices of the subcrystal $\Kup[a]$ and the coordinate
of $v'$ in $\Piup[a]$ is obtained from the one of $v_f$ by increasing its
$k$-th entry by 1. According to Proposition~\ref{pr:fund_string} (with $n$
replaced by $n-1$), $v'$ is obtained from $v_f$ by applying the operator string
  $$
  S_{n-1,k}=w_{n-1,k,n-k}\cdots w_{n-1,k,1},
  $$
where $w_{n-1,k,j}=F_j\cdots F_{j+k-1}$ (cf.~\refeq{string}). In light of this,
we have to show that
  \begin{numitem1}
when (the sequence of moves corresponding to) $S_{n-1,k}$ is applied to $f$,
the resulting feasible function is exactly $f_{a,\Delta'}$.
  \label{eq:applW}
  \end{numitem1}

For convenience $m$-th term $F_{j+m-1}$ in the substring $w_{n-1,k,j}$ will be
denoted by $\phi(j,m)$, $m=1,\ldots,k$. So
$w_{n-1,k,j}=\phi(j,1)\phi(j,2)\ldots \phi(j,k)$.

We distinguish between two cases: $\Delta_k\ge 0$ and $\Delta_k<0$.

 \medskip
{\bf Case 1}: $\Delta_k\ge 0$. An essential fact is that the number $k(n-k)$ of
operators in $S_{n-1,k}$ is equal to the number of nodes in the left rectangle
$L^k$ of $G^k$. Moreover, the level of each node $x_m(j)$ of $L^k$ is equal to
the ``color'' of the operator $\phi(j,m)$ (indeed, $x_m(j)=v^k_{j+m-1}(j)$ and
$\phi(j,m)=F_{j+m-1}$).

Let $f^{j,m}$ denote the current function on $V(G)$ just before the application
of $\phi(j,m)$ (when the process starts with $f=f_{a,\Delta}$). We assert that
  \begin{numitem1}
for each $m$, the application of $\phi(j,m)$ to $f^{j,m}$ increases the value
at the node $x^k_m(j)$ by 1,
  \label{eq:phijm}
 \end{numitem1}
whence~\refeq{applW} will immediately follow.

In order to show~\refeq{phijm}, we first examine tight edges and the slacks
$\eps(v)$ of the nodes $v$ in levels $<n$ for the initial function $f$. One can
see from~\refeq{valuesZ} that
  \begin{numitem1}
for $k'=1,\ldots,n$, each node $v$ of the subgraph $G^{k'}$ has at least one
entering edge (i.e. $\eSW(v)$ or $\eNW(v)$) which is $f$-tight, except,
possibly, for the nodes $v^{k'}_{k'}(1)$, $v^{k'}_{k'-1}(1)$,
$v^{k'}_{n}(n-k'+1)$, $v^{k'}_{n-1}(n-k'+1)$ (indicated by stars in
Fig.~\ref{fig:partition}).
  \label{eq:tightenter}
  \end{numitem1}

 \noindent
{\bf Claim.} {\em For $k'=1,\ldots,n$ and a node $v$ of \,$G^{k'}$ at a level
$<n$,

{\rm(a)} if $v\ne v^{k'}_{k'}(1),v^{k'}_{k'-1}(1)$, then $\eps(v)=0$;

{\rm(b)} if $v= v^{k'}_{k'}(1)$, then $\eps(v)=c_{k'}-a_{k'}-\Delta^+_{k'}\ge
0$;

{\rm(c)} if $v=v^{k'}_{k'-1}(1)$, then $\eps(v)=-\Delta^-_{k'}\ge 0$. }

 \medskip
 \noindent
{\bf Proof of Claim.} Let $v=v^{k'}_i(j)$ and $i<n$. By~\refeq{epsv}, the slack
$\eps(v)$ is equal to $f(w)+f(z)-f(u)-f(v)$, where $w:=v^{k'}_{i-1}(j-1)$,
$z:=v^{k'}_{i+1}(j)$, $u:=v^{k'}_{i}(j-1)$ (these vertices belong to the
extended graph $\bar G^{k'}$). We consider possible cases and
use~\refeq{valuesZ}.

(i) If $w,z,u$ are in $G^{k'}$, then $\partial f(w,v)=\partial f(u,z)$.

(ii) If both $v,w$ are in the piece $Z^{k'}_1$ of $G^{k'}$, then $f(w)=f(v)$
and $f(u)=f(z)=c_{k'}$.

(iii) If $j=1$ and $i\le k'-2$, then $f(v)=f(z)$ and $f(u)=f(w)=c_{k'}$. So in
these cases we have $\eps(v)=0$, yielding~(a).

(iv) Let $v=v^{k'}_{k'}(1)$. Then $f(v)=a_{k'}+\Delta^+_{k'}$ and
$f(u)=f(w)=f(z)=c_{k'}$. This gives $\eps(v)=c_{k'}-a_{k'}-\Delta^+_{k'}$,
yielding~(b).

(v) Let $v=v^{k'}_{k'-1}(1)$. Then $f(v)=a_{k'}+\Delta^+_{k'}+\Delta^-_{k'-1}$,
$f(z)=a_{k'}+\Delta^+_{k'}$ and $f(u)=f(w)=c_{k'}$. This gives
$\eps(v)=-\Delta^-_{k'-1}$, yielding~(c). \hfill \qed

 \medskip
This Claim and the relations $\Delta^-_k=0$ and $\Delta_k<c_k-a_k$ enable us to
estimate the total slacks $\eps_i(j)$ for $f$ at the multinodes $V_i(j)$ with
$i<n$:
  \begin{numitem1}
(i) the edge $\eSW(v^{k+1}_k(1))$ is $f$-tight, $\eps(v^k_k(1))>0$, and
$\eps(v)=0$ for the other \\
 \hphantom{(i)} nodes $v$ in $V_k(1)$; so $\eps_k(1)>0$;

(ii) if $i\ne k,n$, then $\eps(v_i^i(1)),\eps(v_i^{i+1}(1))\ge 0$ and
$\eps(v)=0$ for the other nodes \\
 \hphantom{(ii)} $v$ in $V_i(1)$; so $\eps_i(1)\ge 0$;

(iii) if $i\ne n$ and $j>1$, then $\eps(v)=0$ for all nodes $v$ in $V_i(j)$; so
$\eps_i(j)=0$.
  \label{eq:incase1}
  \end{numitem1}

Now we are ready to prove~\refeq{phijm}. When dealing with a current function
$f^{j,m}$ and seeking for the node at level $j+m-1$ where the operator
$\phi(j,m)$ should act to increase $f^{j,m}$, we can immediately exclude from
consideration any node $v$ that has at least one \emph{tight} entering edge (in
view of the monotonicity condition~\refeq{feas}(i)).

Due to~\refeq{tightenter} and~\refeq{incase1}(i), for the initial function
$f=f^{1,k}$, there is only one node in level $k$ that has no tight entering
edge, namely, $v^k_k(1)$. So, at the first step of the process, the first
operator $\phi(1,k)$ of $S_{n-1,k}$ acts just at $v^k_k(1)$, as required
in~\refeq{phijm}.

Next consider a step with $f':=f^{j,m}$ for $(j,m)\ne (1,k)$, assuming validity
of~\refeq{phijm} on the previous steps.

(A) Let $j=1$ (and $m<k$). For $v:=v^k_m(1)$ and $z:=v^k_{m+1}(1)$, we have
$f'(v)=f(v)\le f(z)=f'(z)-1$. So the unique edge $e=(z,v)$ entering $v$ is not
$f'$-tight. By~\refeq{tightenter}, there are at most two other nodes in level
$m$ that may have no tight entering edges for $f$ (and therefore, for $f'$),
namely, $v^m_m(1)$ and $v^{m+1}_m(1)$. Then $\phi(1,m)$ must act at $v$, as
required in~\refeq{phijm} (since the non-tightness of the SW-edge $e$ of $v$
implies that none of the nodes $v^{m'}_m(1)$ in $V_m(1)$ preceding $v$ (i.e.
with $m'<k$) can be the switch-node).

 \smallskip
(B) Let $j>1$. Comparing $f'$ with $f$ in the node $v:=x^k_m(j)=v^k_{j+m-1}(j)$
and its adjacent nodes, we observe that $v$ has no $f'$-tight entering edge and
that $\eps_{f'}(v)>0$. Also for any other node $v'$ in level $j+m-1$, one can
see that if $v'$ has a tight entering edge for $f$, then so does for $f'$, and
that $\eps_f(v')\ge \eps_{f'}(v')\ge 0$. Using this,
properties~\refeq{tightenter},\,\refeq{incase1}(iii), and
condition~\refeq{activej}, one can conclude that the total and reduced slacks
for $f'$ at the multinode $V':=V_{j+m-1}(j)$ are positive, that $V'$ is the
active multinode for $f'$ in level $j+m-1$, and that $\phi(j,m)$ can be applied
only at $v$, yielding~\refeq{phijm} again.

Thus, \refeq{applW} is valid in Case~1.

 \medskip
{\bf Case 2}: $\Delta_k<0$. We assert that in this case the string $S_{n-1,k}$
acts within the right rectangle $R^{k+1}$ of the subgraph $G^{k+1}$ (note that
$R^{k+1}$ is of size $k\times(n-k)$). More precisely,
  \begin{numitem1}
each operator $\phi(j,m)$ modifies the current function by increasing its value
at the node $x^{k+1}_m(j)$ by 1.
  \label{eq:phijm2}
  \end{numitem1}
Then the resulting function in the process is just $f_{a,\Delta'}$ (in view of
$(\Delta')^-_k=\Delta^-_k+1$), yielding~\refeq{applW}.

To show~\refeq{phijm2}, we argue as in the previous case and
use~\refeq{tightenter} and the above Claim. Since $\Delta_k<0$, part~(i)
in~\refeq{incase1} for the initial function $f$ is modified as:
  \begin{numitem1}
for $j=1,\ldots,n-k$, the SW-edge of each node
$x^{k+1}_k(j)=v^{k+1}_{j+k-1}(j)$ is not $f$-tight, $\eps(v^{k+1}_k(1))>0$,
~$\eps(v^k_k(1))\ge 0$, and $\eps(v)=0$ for the other nodes $v$ in $V_k(1)$; so
$\eps_k(1)>0$,
  \label{eq:incase2}
  \end{numitem1}
while properties~(ii) and~(iii) preserve.

By~\refeq{tightenter} and~\refeq{incase2}, there are only two nodes in level
$k$ that have no $f$-tight entering edges, namely, $v^k_k(1)$ and
$v^{k+1}_k(1)$. Also $e=\eSW(v^{k+1}_k(1))$ is not tight. So, at the first
step, $\phi(1,k)$ must act at $v^{k+1}_k(1)$, as required in~\refeq{phijm2}
(since the non-tightness of $e$ implies that the node $v^k_k(1)$ preceding
$v^{k+1}_k(1)$ cannot be the switch-node in $V_k(1)$).

The fact that $\phi(1,m)$ with $m<k$ acts at $v^{k+1}_m(1)$ is shown by arguing
as in~(A) above. And for $j>1$, to show that $\phi(j,m)=F_{j+m-1}$ acts at
$x^{k+1}_m(j)=v^{k+1}_{j+m-1}(j)$, we argue as in~(B) above. Here, when $m=k$,
we also use the fact that the edge $\eSW(x^{k+1}_{k}(j))$ is not $f$-tight
(by~\refeq{incase2}), whence both edges entering $x^{k+1}_k(j)$ are not tight
for the current function. So~\refeq{phijm2} is always valid.

Thus, we have the desired property~\refeq{applW} in both cases~1 and~2, and
statement~(ii) in Lemma~\ref{lm:faDelta} follows. \hfill \qed\qed

This completes the proof of relation~\refeq{bi}, yielding
Theorem~\ref{tm:two_loci}.

 \section{Illustrations} \label{sec:illustr}

In this concluding section we give two illustrations to the above assembling
construction for A-crystals. The first one refines the interrelation between
upper and lower subcrystals in an arbitrary $A_2$-crystal; this can be compared
with the explicit construction (the so-called ``sail model'') for
$A_2$-crystals in~\cite{A2}. The second one visualizes the subcrystals
structure for one instance of $A_3$-crystals, namely, $K(1,1,1)$.

\subsection{$A_2$-crystals} \label{ssec:A2}

The subcrystals structure becomes simpler when we deal with an $A_2$-crystal
$K=K(c_1,c_2)$. In this case the roles of upper, lower, and middle subcrystals
are played by 1-paths, 2-paths, and vertices of $K$, respectively, where by an
$i$-path we mean a maximal path of color $i$.

Consider an upper subcrystal in $K$. This is a 1-path $P=(v_0,v_1,\ldots,v_p)$
containing exactly one principal vertex $\prv[a]$ of $K$ (the heart of $P$);
here $v_i$ stands for $i$-th vertex in $P$, ~$a=(a_1,a_2)\in\Zset_+^2$, and
$a\le c$. Let $\prv[a]=v_h$. Formulas~\refeq{par_up} and~\refeq{heart_up} give
  \begin{equation} \label{eq:1path}
  |P|=p=c_1-a_1+a_2 \qquad\mbox{and} \qquad h=a_2.
  \end{equation}
Fix a vertex $v=v_i$ of $P$. It belongs to some 2-path (lower subcrystal)
$Q=(u_1,u_2,\ldots,u_q)$. Let $v=u_j$ and let $\prv[b]=u_{\bar h}$ be the
principal vertex of $K$ occurring in $Q$ (the heart of $Q$). The vertex $v$
forms a middle subcrystal of $K$; its deviations from the heart of $P$ and from
the heart of $Q$ are equal to $i-h=:\delta$ and $j-\bar h=:\bar\delta$,
respectively. By~\refeq{nablai} in Theorem~\ref{tm:two_div}, we have
$\bar\delta=-\delta$. Then we can compute the coordinates $b$ by use
of~\refeq{bi} and, further, apply~\refeq{par_low} and~\refeq{heart_low} to
compute the length of $Q$ and the locus of its heart. This gives:
  \begin{numitem1} \label{eq:2path}
 \begin{itemize}
\item[(i)] if $\delta\ge 0$ (i.e. $a_2\le i\le c_1-a_1+a_2$), then
$b_1=a_1+\delta=a_1+i-a_2$, ~$b_2=a_2$, ~$|Q|=c_2-b_2+b_1=c_2-2a_2+a_1+i$, and
$|Q|-\bar h=|Q|-b_1=c_2-a_2$;
  \item[(ii)] if $\delta\le 0$ (i.e. $0\le i\le a_2$), then $b_1=a_1$,
~$b_2=a_2+\delta=a_2+(i-a_2)=i$, ~$|Q|=c_2-b_2+b_1=c_2-i+a_1$, and $\bar
h=b_1=a_1$.
  \end{itemize}
  \end{numitem1}

Using~\refeq{1path} and~\refeq{2path}, one can enumerate the sets of 1-paths
and 2-paths and properly intersect corresponding pairs, obtaining the
$A_2$-crystal $K(c)$. It is rather routine to check that the resulting graph
coincides with the one generated by the \emph{sail model} from~\cite{A2}. Next
we outline that construction.

Given $c\in\Zset_+^2$, the $A_2$-crystal $K(c)$ is produced from two particular
two-colored graphs $R$ and $L$, called the {\em right sail} of size $c_1$ and
the {\em left sail} of size $c_2$, respectively. The vertices of $R$ correspond
to the vectors $(i,j)\in\Zset^2$ such that $0\le j\le i\le c_1$, and the
vertices of $L$ to the vectors $(i,j)\in\Zset^2$ such that $0\le i\le j\le
c_2$. In both $R$ and $L$, the edges of color 1 are all possible pairs of the
form $((i,j),(i+1,j))$, and the edges of color 2 are all possible pairs of the
form $((i,j),(i,j+1))$. (Observe that both $R$ and $L$ satisfy axioms
(A1)-(A4), $R$ is isomorphic to $K(c_1,0)$, ~$L$ is isomorphic to $K(0,c_2)$,
and their critical vertices are the ``diagonal vertices'' $(i,i)$.)

In order to produce $K(c)$, take $c_2$ disjoint copies $R_1,\ldots,R_{c_2}$ of
$R$ and $c_1$ disjoint copies $L_1,\ldots,L_{c_1}$ of $L$, referring to $R_j$
as $j$-th right sail, and to $L_i$ as $i$-th left sail. Let $D(R_j)$ and
$D(L_i)$ denote the sets of diagonal vertices in $R_j$ and $L_i$, respectively.
For all $i=1,\ldots,c_1$ and $j=1,\ldots,c_2$, we identify the diagonal
vertices $(i,i)\in D(R_j)$ and $(j,j)\in D(L_i)$. The resulting graph is just
the desired $K(c)$. The edge colors of $K(c)$ are inherited from $L$ and $R$.
One checks that $K(c)$ has $(c_1+1)\times (c_2+1)$ critical vertices; they
coincide with the diagonal vertices of the sails. The principal lattice of
$K(c)$ is just constituted by the critical vertices.

The case $(c_1,c_2)=(1,2)$ is drawn in the picture; here the critical
(principal) vertices are indicated by big circles, 1-edges by horizontal
arrows, and 2-edges by vertical arrows.
 \begin{center}
  \unitlength=1mm
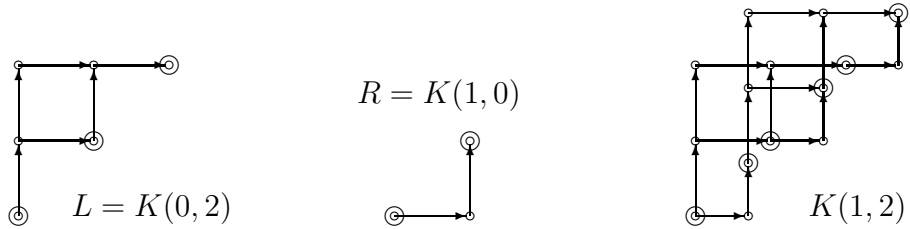
  \begin{picture}(125,35)
\put(5,5){\circle{1.0}}           
\put(5,15){\circle{1.0}} \put(5,25){\circle{1.0}} \put(15,15){\circle{1.0}}
\put(15,25){\circle{1.0}} \put(25,25){\circle{1.0}} \put(5,5){\circle{2.5}}
\put(15,15){\circle{2.5}} \put(25,25){\circle{2.5}}
\put(5,15){\vector(1,0){9.5}} \put(5,25){\vector(1,0){9.5}}
\put(15,25){\vector(1,0){9.5}} \put(5,5){\vector(0,1){9.5}}
\put(5,15){\vector(0,1){9.5}} \put(15,15){\vector(0,1){9.5}}
\put(12,5){$L=K(0,2)$}
\put(55,5){\circle{1.0}} \put(65,5){\circle{1.0}} \put(65,15){\circle{1.0}}
\put(55,5){\circle{2.5}} \put(65,15){\circle{2.5}}
\put(55,5){\vector(1,0){9.5}} \put(65,5){\vector(0,1){9.5}}
 \put(50,20){$R=K(1,0)$}
\put(95,5){\circle{1.0}} \put(95,15){\circle{1.0}} \put(95,25){\circle{1.0}}
\put(105,15){\circle{1.0}} \put(105,25){\circle{1.0}}
\put(115,25){\circle{1.0}} \put(95,5){\circle{2.5}} \put(105,15){\circle{2.5}}
\put(115,25){\circle{2.5}} \put(95,15){\vector(1,0){9.5}}
\put(95,25){\vector(1,0){9.5}} \put(105,25){\vector(1,0){9.5}}
\put(95,5){\vector(0,1){9.5}} \put(95,15){\vector(0,1){9.5}}
\put(105,15){\vector(0,1){9.5}}
\put(102,12){\circle{1.0}} \put(102,22){\circle{1.0}}
\put(102,32){\circle{1.0}} \put(112,22){\circle{1.0}}
\put(112,32){\circle{1.0}} \put(122,32){\circle{1.0}}
\put(102,12){\circle{2.5}} \put(112,22){\circle{2.5}}
\put(122,32){\circle{2.5}} \put(102,22){\vector(1,0){9.5}}
\put(102,32){\vector(1,0){9.5}} \put(112,32){\vector(1,0){9.5}}
\put(102,12){\vector(0,1){9.5}} \put(102,22){\vector(0,1){9.5}}
\put(112,22){\vector(0,1){9.5}}
\put(102,5){\circle{1.0}} \put(95,5){\vector(1,0){6.5}}
\put(102,5){\vector(0,1){6.5}} \put(112,15){\circle{1.0}}
\put(105,15){\vector(1,0){6.5}} \put(112,15){\vector(0,1){6.5}}
\put(122,25){\circle{1.0}} \put(115,25){\vector(1,0){6.5}}
\put(122,25){\vector(0,1){6.5}}
 \put(110,5){$K(1,2)$}
  \end{picture}
 \end{center}

In particular, the sail model shows that the numbers of edges of each color in
an $A_2$-crystal are the same. This implies a similar property for any
$A_n$-crystal.

\subsection{$A_3$-crystal $K(1,1,1)$} \label{ssec:K111}

Next we illustrate the $A_3$-crystal $K=K(1,1,1)$. It has 64 vertices and 102
edges, is rather puzzling, and drawing it in full would take too much space;
for this reason, we show it by fragments, namely, by demonstrating all of its
upper and lower subcrystals. We abbreviate notation $\prv[(i,j,k)]$ for
principal vertices to $(i,j,k)$ for short. So the principal lattice consists of
eight vertices $(0,0,0),\ldots,(1,1,1)$, as drawn in the picture (where the
arrows indicate moves by principal operator strings $S_{3,k}$ as
in~\refeq{string}):
  \begin{center}
  \unitlength=1mm
  \begin{picture}(80,42)
    \put(5,5){\begin{picture}(80,35)
  \put(0,0){\circle*{1.2}}
  \put(30,0){\circle*{1.2}}
  \put(0,20){\circle*{1.2}}
  \put(10,10){\circle*{1.2}}
  \put(40,10){\circle*{1.2}}
  \put(10,30){\circle*{1.2}}
  \put(30,20){\circle*{1.2}}
  \put(40,30){\circle*{1.2}}
  \put(0,0){\vector(1,0){29}}
  \put(10,10){\vector(1,0){29}}
  \put(0,20){\vector(1,0){29}}
  \put(10,30){\vector(1,0){29}}
  \put(0,0){\vector(0,1){19.3}}
  \put(30,0){\vector(0,1){19.3}}
  \put(10,10){\vector(0,1){19.3}}
  \put(40,10){\vector(0,1){19.3}}
  \put(0,0){\vector(1,1){9.5}}
  \put(30,0){\vector(1,1){9.5}}
  \put(0,20){\vector(1,1){9.5}}
  \put(30,20){\vector(1,1){9.5}}
  \put(-14,-4){$s$=(0,0,0)}
  \put(28,-4){(1,0,0)}
  \put(-12,21){(0,1,0)}
  \put(11,11){(0,0,1)}
  \put(19,21){(1,1,0)}
  \put(41,8){(1,0,1)}
  \put(4,31.5){(0,1,1)}
  \put(33,31.5){(1,1,1)=$t$}
  \put(10,-4){$S_{3,1}$}
  \put(-7,9){$S_{3,2}$}
  \put(6,3){$S_{3,3}$}
  \put(60,20){$S_{3,1}=F_3F_2F_1$}
  \put(60,13){$S_{3,2}=F_2F_3F_1F_2$}
  \put(60,6){$S_{3,3}=F_1F_2F_3$}
   \end{picture}}
 \end{picture}
 \end{center}

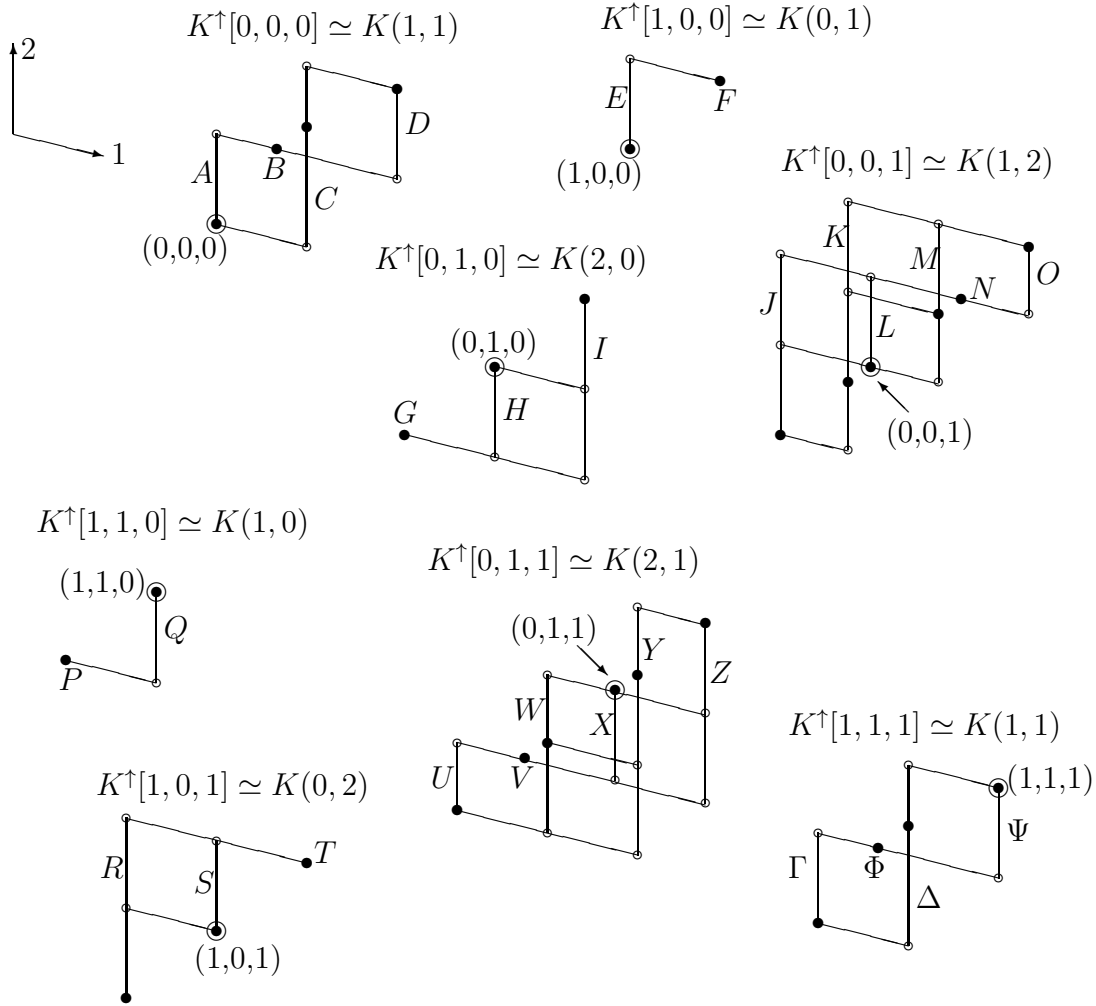
\begin{figure}[hbt]                  
  \begin{center}
  \unitlength=1mm
  \begin{picture}(145,137)
   \put(3,105){\begin{picture}(20,30)
  \put(0,10){\vector(4,-1){12}}
  \put(0,10){\vector(0,1){12}}
  \put(13,6){1}
  \put(1,20){2}
   \end{picture}}
   \put(30,103){\begin{picture}(40,30)
  \put(0,0){\circle*{1.5}}
  \put(12,-3){\circle{1.2}}
  \put(0,12){\circle{1.2}}
  \put(24,6){\circle{1.2}}
  \put(8,10){\circle*{1.5}}
  \put(12,21){\circle{1.2}}
  \put(12,13){\circle*{1.5}}
  \put(24,18){\circle*{1.5}}
  \put(0,0){\circle{2.5}}
  \put(0,0){\line(4,-1){12}}
  \put(0,12){\line(4,-1){24}}
  \put(12,21){\line(4,-1){12}}
  \put(0,0){\line(0,1){12}}
  \put(12,-3){\line(0,1){24}}
  \put(24,6){\line(0,1){12}}
  \put(-10,-4.5){(0,0,0)}
  \put(-3.5,5.5){$A$}
  \put(6,6){$B$}
  \put(13,2){$C$}
  \put(25,12){$D$}
  \put(-4,25){$\Kup[0,0,0]\simeq K(1,1)$}
   \end{picture}}
   \put(85,113){\begin{picture}(30,25)
  \put(0,0){\circle*{1.5}}
  \put(0,12){\circle{1.2}}
  \put(12,9){\circle*{1.5}}
  \put(0,0){\circle{2.5}}
  \put(0,12){\line(4,-1){12}}
  \put(0,0){\line(0,1){12}}
  \put(-10,-4.5){(1,0,0)}
  \put(-3.5,5.5){$E$}
  \put(11,5){$F$}
  \put(-4,16){$\Kup[1,0,0]\simeq K(0,1)$}
   \end{picture}}
   \put(55,75){\begin{picture}(40,30)
  \put(0,0){\circle*{1.5}}
  \put(12,-3){\circle{1.2}}
  \put(24,-6){\circle{1.2}}
  \put(24,6){\circle{1.2}}
  \put(12,9){\circle*{1.5}}
  \put(24,18){\circle*{1.5}}
  \put(12,9){\circle{2.5}}
  \put(0,0){\line(4,-1){24}}
  \put(12,9){\line(4,-1){12}}
  \put(12,-3){\line(0,1){12}}
  \put(24,-6){\line(0,1){24}}
  \put(6,11){(0,1,0)}
  \put(-1.5,1.5){$G$}
  \put(13,2){$H$}
  \put(25,10){$I$}
  \put(-4,22){$\Kup[0,1,0]\simeq K(2,0)$}
   \end{picture}}

   \put(105,75){\begin{picture}(45,45)
  \put(0,0){\circle*{1.5}}
  \put(0,12){\circle{1.2}}
  \put(0,24){\circle{1.2}}
  \put(12,9){\circle*{1.5}}
  \put(12,21){\circle{1.2}}
  \put(24,18){\circle*{1.5}}
  \put(12,9){\circle{2.5}}
  \put(0,12){\line(4,-1){21}}
  \put(0,24){\line(4,-1){33}}
  \put(0,0){\line(0,1){24}}
  \put(12,9){\line(0,1){12}}
  \put(9,7){\circle*{1.5}}
  \put(9,19){\circle{1.2}}
  \put(9,31){\circle{1.2}}
  \put(21,16){\circle*{1.5}}
  \put(21,28){\circle{1.2}}
  \put(33,25){\circle*{1.5}}
  \put(9,19){\line(4,-1){12}}
  \put(9,31){\line(4,-1){24}}
  \put(9,-2){\line(0,1){33}}
  \put(21,7){\line(0,1){21}}
  \put(9,-2){\circle{1.2}}
  \put(21,7){\circle{1.2}}
  \put(33,16){\circle{1.2}}
  \put(0,0){\line(4,-1){9}}
  \put(33,16){\line(0,1){9}}
  \put(14,-1){(0,0,1)}
 \put(17,3){\vector(-1,1){4}}
  \put(-3,16){$J$}
  \put(5.5,25){$K$}
  \put(12.5,13){$L$}
  \put(17,22){$M$}
  \put(25,18){$N$}
  \put(34,20){$O$}
  \put(0,35){$\Kup[0,0,1]\simeq K(1,2)$}
   \end{picture}}
   \put(10,45){\begin{picture}(30,25)
  \put(0,0){\circle*{1.5}}
  \put(12,-3){\circle{1.2}}
  \put(12,9){\circle*{1.5}}
  \put(12,9){\circle{2.5}}
  \put(0,0){\line(4,-1){12}}
  \put(12,-3){\line(0,1){12}}
  \put(-1,9){(1,1,0)}
  \put(-1,-4){$P$}
  \put(13,3){$Q$}
  \put(-4,17){$\Kup[1,1,0]\simeq K(1,0)$}
   \end{picture}}

   \put(18,0){\begin{picture}(40,30)
  \put(0,0){\circle*{1.5}}
  \put(0,12){\circle{1.2}}
  \put(0,24){\circle{1.2}}
  \put(12,9){\circle*{1.5}}
  \put(12,21){\circle{1.2}}
  \put(24,18){\circle*{1.5}}
  \put(12,9){\circle{2.5}}
  \put(0,12){\line(4,-1){12}}
  \put(0,24){\line(4,-1){24}}
  \put(0,0){\line(0,1){24}}
  \put(12,9){\line(0,1){12}}
  \put(9,4){(1,0,1)}
  \put(-3.5,16){$R$}
  \put(9,14){$S$}
  \put(25,18){$T$}
  \put(-4,27){$\Kup[1,0,1]\simeq K(0,2)$}
   \end{picture}}

   \put(62,25){\begin{picture}(50,50)
  \put(0,0){\circle*{1.5}}
  \put(12,-3){\circle{1.2}}
  \put(24,-6){\circle{1.2}}
  \put(24,6){\circle{1.2}}
  \put(12,9){\circle*{1.5}}
  \put(24,18){\circle*{1.5}}
  \put(0,0){\line(4,-1){24}}
  \put(12,9){\line(4,-1){12}}
  \put(12,-3){\line(0,1){21}}
  \put(24,-6){\line(0,1){33}}
  \put(9,7){\circle*{1.5}}
  \put(21,4){\circle{1.2}}
  \put(33,1){\circle{1.2}}
  \put(33,13){\circle{1.2}}
  \put(21,16){\circle*{1.5}}
  \put(33,25){\circle*{1.5}}
  \put(21,16){\circle{2.5}}
  \put(0,9){\line(4,-1){33}}
  \put(12,18){\line(4,-1){21}}
  \put(21,4){\line(0,1){12}}
  \put(33,1){\line(0,1){24}}
  \put(0,9){\circle{1.2}}
  \put(12,18){\circle{1.2}}
  \put(24,27){\circle{1.2}}
  \put(24,27){\line(4,-1){9}}
  \put(0,0){\line(0,1){9}}
  \put(7,23){(0,1,1)}
 \put(16,22){\vector(1,-1){4}}
  \put(-3.5,3){$U$}
  \put(7,3){$V$}
  \put(7.5,12){$W$}
  \put(17.5,10){$X$}
  \put(24.5,20){$Y$}
  \put(33.5,17){$Z$}
  \put(-4,32){$\Kup[0,1,1]\simeq K(2,1)$}
   \end{picture}}
   \put(110,10){\begin{picture}(40,30)
  \put(0,0){\circle*{1.5}}
  \put(12,-3){\circle{1.2}}
  \put(0,12){\circle{1.2}}
  \put(24,6){\circle{1.2}}
  \put(8,10){\circle*{1.5}}
  \put(12,21){\circle{1.2}}
  \put(12,13){\circle*{1.5}}
  \put(24,18){\circle*{1.5}}
  \put(24,18){\circle{2.5}}
  \put(0,0){\line(4,-1){12}}
  \put(0,12){\line(4,-1){24}}
  \put(12,21){\line(4,-1){12}}
  \put(0,0){\line(0,1){12}}
  \put(12,-3){\line(0,1){24}}
  \put(24,6){\line(0,1){12}}
  \put(25,18){(1,1,1)}
  \put(-4,6){$\Gamma$}
  \put(6,6){$\Phi$}
  \put(13,2){$\Delta$}
  \put(25,11){$\Psi$}
  \put(-4,25){$\Kup[1,1,1]\simeq K(1,1)$}
   \end{picture}}
 \end{picture}
 \end{center}
  \caption{The upper subcrystals in $K(1,1,1)$}
  \label{fig:upper111}
  \end{figure}

\begin{figure}[hbt]                  
  \begin{center}
  \unitlength=1mm
  \begin{picture}(145,148)
   \put(3,110){\begin{picture}(20,30)
  \put(0,10){\vector(3,1){9}}
  \put(0,10){\vector(0,1){12}}
  \put(1,20){2}
  \put(9.5,13){3}
   \end{picture}}
   \put(30,105){\begin{picture}(40,40)
  \put(0,0){\circle*{1.5}}
  \put(12,4){\circle{1.2}}
  \put(0,12){\circle{1.2}}
  \put(24,20){\circle{1.2}}
  \put(8,14.7){\circle*{1.5}}
  \put(12,28){\circle{1.2}}
  \put(12,20){\circle*{1.5}}
  \put(24,32){\circle*{1.5}}
  \put(0,0){\circle{2.5}}
  \put(0,0){\line(3,1){12}}
  \put(0,12){\line(3,1){24}}
  \put(12,28){\line(3,1){12}}
  \put(0,0){\line(0,1){12}}
  \put(12,4){\line(0,1){24}}
  \put(24,20){\line(0,1){12}}
  \put(-10,-4.5){(0,0,0)}
  \put(-3.5,5.5){$A$}
  \put(6,10){$G$}
  \put(13,9){$J$}
  \put(25,24.5){$U$}
  \put(-4,35){$\Klow[0,0,0]\simeq K(1,1)$}
   \end{picture}}
   \put(85,120){\begin{picture}(30,25)
  \put(0,0){\circle*{1.5}}
  \put(0,12){\circle{1.2}}
  \put(12,16){\circle*{1.5}}
  \put(0,0){\circle{2.5}}
  \put(0,12){\line(3,1){12}}
  \put(0,0){\line(0,1){12}}
  \put(-10,-4.5){(0,0,1)}
  \put(-3.5,5.5){$L$}
  \put(12,11.5){$V$}
  \put(-4,20){$\Klow[0,0,1]\simeq K(0,1)$}
   \end{picture}}
   \put(53,70){\begin{picture}(40,35)
  \put(0,0){\circle*{1.5}}
  \put(12,4){\circle{1.2}}
  \put(24,8){\circle{1.2}}
  \put(24,20){\circle{1.2}}
  \put(12,16){\circle*{1.5}}
  \put(24,32){\circle*{1.5}}
  \put(12,16){\circle{2.5}}
  \put(0,0){\line(3,1){24}}
  \put(12,16){\line(3,1){12}}
  \put(12,4){\line(0,1){12}}
  \put(24,8){\line(0,1){24}}
  \put(3,18){(0,1,0)}
  \put(-1.5,1.5){$B$}
  \put(12.5,10){$H$}
  \put(25,24){$W$}
  \put(-4,35){$\Klow[0,1,0]\simeq K(2,0)$}
   \end{picture}}

   \put(103,75){\begin{picture}(45,50)
  \put(0,0){\circle*{1.5}}
  \put(0,12){\circle{1.2}}
  \put(0,24){\circle{1.2}}
  \put(12,16){\circle*{1.5}}
  \put(12,28){\circle{1.2}}
  \put(24,32){\circle*{1.5}}
  \put(12,16){\circle{2.5}}
  \put(0,12){\line(3,1){21}}
  \put(0,24){\line(3,1){33}}
  \put(0,0){\line(0,1){24}}
  \put(12,16){\line(0,1){12}}
  \put(9,11.7){\circle*{1.5}}
  \put(9,23.7){\circle{1.2}}
  \put(9,35.7){\circle{1.2}}
  \put(21,28.2){\circle*{1.5}}
  \put(21,40.2){\circle{1.2}}
  \put(33,44.2){\circle*{1.5}}
  \put(9,24.2){\line(3,1){12}}
  \put(9,36.2){\line(3,1){24}}
  \put(9,3.2){\line(0,1){33}}
  \put(21,19.2){\line(0,1){21}}
  \put(9,3.2){\circle{1.2}}
  \put(21,19.2){\circle{1.2}}
  \put(33,35.2){\circle{1.2}}
  \put(0,0){\line(3,1){9}}
  \put(33,35.2){\line(0,1){9}}
  \put(15,7){(1,0,0)}
 \put(17.5,10.5){\vector(-1,1){4}}
  \put(-3.5,16){$C$}
  \put(5,30){$K$}
  \put(12.5,20){$E$}
  \put(17.5,34){$R$}
  \put(24.5,29){$P$}
  \put(34,39){$\Gamma$}
  \put(0,48){$\Klow[1,0,0]\simeq K(1,2)$}
   \end{picture}}
   \put(5,45){\begin{picture}(30,27)
  \put(0,0){\circle*{1.5}}
  \put(12,4){\circle{1.2}}
  \put(12,16){\circle*{1.5}}
  \put(12,16){\circle{2.5}}
  \put(0,0){\line(3,1){12}}
  \put(12,4){\line(0,1){12}}
  \put(-1,16){(0,1,1)}
  \put(-1,-4){$N$}
  \put(13,10){$X$}
  \put(-4,23){$\Klow[0,1,1]\simeq K(1,0)$}
   \end{picture}}

   \put(15,0){\begin{picture}(40,30)
  \put(0,0){\circle*{1.5}}
  \put(0,12){\circle{1.2}}
  \put(0,24){\circle{1.2}}
  \put(12,16){\circle*{1.5}}
  \put(12,28){\circle{1.2}}
  \put(24,32){\circle*{1.5}}
  \put(12,16){\circle{2.5}}
  \put(0,12){\line(3,1){12}}
  \put(0,24){\line(3,1){24}}
  \put(0,0){\line(0,1){24}}
  \put(12,16){\line(0,1){12}}
  \put(12,12){(1,0,1)}
  \put(-4,16){$M$}
  \put(12.5,21){$S$}
  \put(24,28){$\Phi$}
  \put(-4,35){$\Klow[1,0,1]\simeq K(0,2)$}
   \end{picture}}

   \put(62,13){\begin{picture}(50,50)
  \put(0,0){\circle*{1.5}}
  \put(12,4){\circle{1.2}}
  \put(24,8){\circle{1.2}}
  \put(24,20){\circle{1.2}}
  \put(12,16){\circle*{1.5}}
  \put(24,32){\circle*{1.5}}
  \put(0,0){\line(3,1){24}}
  \put(12,16){\line(3,1){12}}
  \put(12,4){\line(0,1){21}}
  \put(24,8){\line(0,1){33}}
  \put(9,12.1){\circle*{1.5}}
  \put(21,16.2){\circle{1.2}}
  \put(33,20.2){\circle{1.2}}
  \put(33,32.2){\circle{1.2}}
  \put(21,28.2){\circle*{1.5}}
  \put(33,44.2){\circle*{1.5}}
  \put(21,28.2){\circle{2.5}}
  \put(0,9){\line(3,1){33}}
  \put(12,25){\line(3,1){21}}
  \put(21,16.2){\line(0,1){12}}
  \put(33,20.2){\line(0,1){24}}
  \put(0,9){\circle{1.2}}
  \put(12,25){\circle{1.2}}
  \put(24,41){\circle{1.2}}
  \put(24,41){\line(3,1){9}}
  \put(0,0){\line(0,1){9}}
  \put(9,29.5){(1,1,0)}
  \put(-3.5,3){$D$}
  \put(7,8){$F$}
  \put(9.5,19){$I$}
  \put(17.5,21.5){$Q$}
  \put(24.5,35){$Y$}
  \put(33.5,38){$\Delta$}
  \put(-4,46){$\Klow[1,1,0]\simeq K(2,1)$}
   \end{picture}}
   \put(110,5){\begin{picture}(40,30)
  \put(0,0){\circle*{1.5}}
  \put(12,4){\circle{1.2}}
  \put(0,12){\circle{1.2}}
  \put(24,20){\circle{1.2}}
  \put(8,14.7){\circle*{1.5}}
  \put(12,28){\circle{1.2}}
  \put(12,20){\circle*{1.5}}
  \put(24,32){\circle*{1.5}}
  \put(24,32){\circle{2.5}}
  \put(0,0){\line(3,1){12}}
  \put(0,12){\line(3,1){24}}
  \put(12,28){\line(3,1){12}}
  \put(0,0){\line(0,1){12}}
  \put(12,4){\line(0,1){24}}
  \put(24,20){\line(0,1){12}}
  \put(25,31){(1,1,1)}
  \put(-3.5,5){$O$}
  \put(7,10){$T$}
  \put(8.5,22){$Z$}
  \put(24.5,24){$\Psi$}
  \put(0,38){$\Klow[1,1,1]\simeq K(1,1)$}
   \end{picture}}
 \end{picture}
 \end{center}
  \caption{The lower subcrystals in $K(1,1,1)$}
  \label{fig:lower111}
  \end{figure}
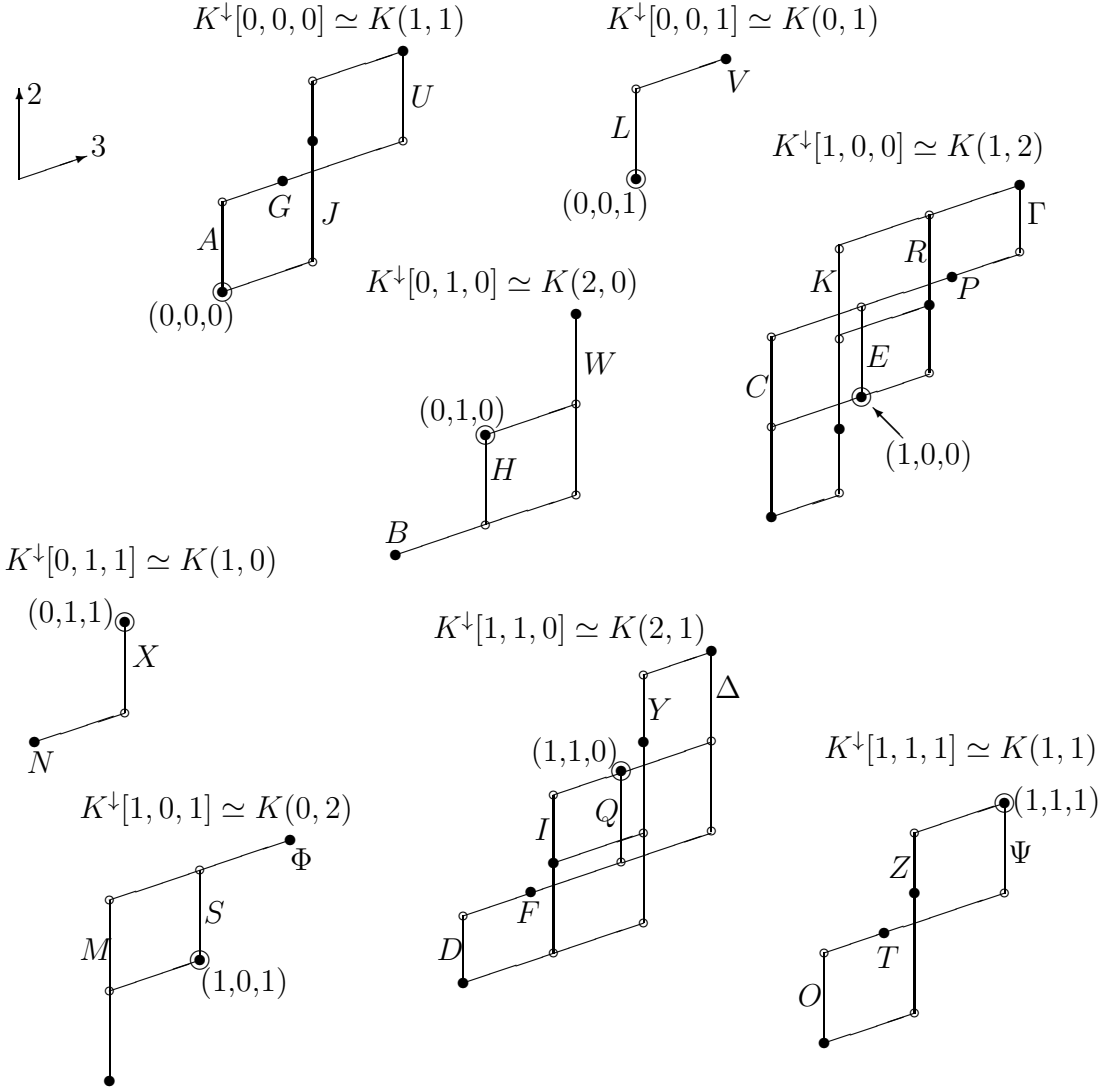

Thus, $K$ has eight upper subcrystals $\Kup[i,j,k]$ and eight lower subcrystals
$\Klow[i,j,k]$ (writing $K^\bullet[i,j,k]$ for $K^\bullet[(i,j,k)]$); they are
drawn in Figures~\ref{fig:upper111} and~\ref{fig:lower111}. Here the directions
of edges of colors 1,2,3 are as indicated in the upper left corner. In each
subcrystal we indicate its critical vertices by black circles, and the unique
principal vertex of $K$ occurring in it (the heart) by a big white circle. $K$
has 30 middle subcrystals (paths of color 2), which are labeled as $A,\ldots,Z,
\Gamma,\Delta,\Phi,\Psi$ (note that $B,F,G,N,P,T,V,\Phi$ consist of single
vertices).

For each upper subcrystal $\Kup[i,j,k]$, its parameter $\parup$ and heart locus
$\heartup$, computed by use of~\refeq{par_up} and~\refeq{heart_up}, are as
follows (where $\tilde K,\tilde s,z$ denote the current subcrystal, its source,
and its heart, respectively):

$\bullet$ for $\Kup[0,0,0]$: $\parup_1=1-0+0=1$, $\parup_2=1-0+0=1$, and
$\heartup_1=\heartup_2=0$ (so $\tilde K$ is isomorphic to $K(1,1)$ and $z$
coincides with $\tilde s$);

$\bullet$ for $\Kup[1,0,0]$: $\parup_1=1-1+0=0$, $\parup_2=1-0+0=1$, and
$\heartup_1=\heartup_2=0$;

$\bullet$ for $\Kup[0,1,0]$: $\parup_1=1-0+1=2$, $\parup_2=1-1+0=0$,
$\heartup_1=1$, and $\heartup_2=0$ (so $\tilde K\simeq K(2,0)$ and $z$ is
located at $S_{2,1}(\tilde s)=F_2F_1(\tilde s)$);

$\bullet$ for $\Kup[0,0,1]$: $\parup_1=1-0+0=1$, $\parup_2=1-0+1=2$,
$\heartup_1=0$, and $\heartup_2=1$ (so $\tilde K\simeq K(1,2)$ and $z$ is
located at $S_{2,2}(\tilde s)=F_1F_2(\tilde s)$);

$\bullet$ for $\Kup[1,1,0]$: $\parup_1=1-1+1=1$, $\parup_2=1-1+0=0$,
$\heartup_1=1$, and $\heartup_2=0$;

$\bullet$ for $\Kup[1,0,1]$: $\parup_1=1-1+0=0$, $\parup_2=1-0+1=2$,
$\heartup_1=0$, and $\heartup_2=1$;

$\bullet$ for $\Kup[0,1,1]$: $\parup_1=1-0+1=2$, $\parup_2=1-1+1=1$, and
$\heartup_1=\heartup_2=1$ (so $\tilde K\simeq K(2,1)$ and $z$ is located at
$S_{2,2}S_{2,1}(\tilde s)=F_1F_2F_2F_1(\tilde s)$);

$\bullet$ for $\Kup[1,1,1]$: $\parup_1=1-1+1=1$, $\parup_2=1-1+1=1$, and
$\heartup_1=\heartup_2=1$.

Since $K(1,1,1)$ is ``symmetric'', so are its upper and lower subcrystals, i.e.
each $\Klow[i,j,k]$ is obtained from $\Kup[k,j,i]$ by replacing color 1 by 3.
In Fig.~\ref{fig:lower111}, when writing $\Klow[i,j,k]\simeq K(\alpha,\beta)$,
the parameters $\alpha,\beta$ concern colors 3 and 2, respectively.

Now the desired $K(1,1,1)$ is assembled by gluing the fragments in
Figs.~\ref{fig:upper111},\ref{fig:lower111} along the 2-paths $A,\ldots,\Psi$.



\end{document}